\documentclass[journal]{IEEEtran}
\usepackage{amsfonts}
\usepackage{bm}
\usepackage{amsmath,amsfonts,amsthm,amssymb}
\usepackage{setspace}
\usepackage{Tabbing}
\usepackage{fancyhdr}
\usepackage{lastpage}
\usepackage{extramarks}
\usepackage{chngpage}
\usepackage{algorithmic} 
\usepackage{soul,color}
\usepackage{epsfig}
\usepackage{graphicx,float,wrapfig,subfigure}
\usepackage{epstopdf}
\usepackage{dsfont}
\usepackage{longtable}
\usepackage{wrapfig}
\usepackage{stfloats}
\usepackage{cite}
\usepackage{graphicx}
\usepackage{array}
\usepackage{multirow}{\tiny }
\usepackage{multicol}
\usepackage{colortbl}
\usepackage{tabularx}
\usepackage{mdwmath}
\usepackage{mdwtab}
\usepackage{color}
\usepackage{verbatim}
\usepackage{amsmath}
\usepackage{tikz}
\usetikzlibrary{calc}
\usepackage{flowchart}
\usetikzlibrary{shapes.geometric, arrows}
\usepackage{overpic}
\usepackage{url}
\usepackage[colorlinks,linkcolor=black,anchorcolor=black,citecolor=black,urlcolor=black]{hyperref} 
\usepackage{breakurl}
\usepackage[flushleft]{threeparttable}
\allowdisplaybreaks

\newcounter{remark}
\newenvironment{remark}[1][]{\refstepcounter{remark}\par
	\textbf{Remark~\theremark #1:} \rmfamily}

\newcounter{proposition}
\newenvironment{proposition}[1][]{\refstepcounter{proposition}\par
	\textbf{Proposition~\theproposition #1:} \rmfamily}

\begin{document}

\title{
	Data-driven Chance-constrained Regulation Capacity Offering for Distributed Energy Resources
	%Under Dispatch Uncertainty
	}

\author{
Hongcai~Zhang,~\IEEEmembership{Student Member,~IEEE,}
Zechun~Hu,~\IEEEmembership{Senior Member,~IEEE,}
Eric~Munsing,~\IEEEmembership{Student Member,~IEEE,}
Scott~J.~Moura,~\IEEEmembership{Member,~IEEE,}
and~Yonghua~Song,~\IEEEmembership{Fellow,~IEEE
	\vspace{-2mm}}

\thanks{
%	This work was supported in part by the National Key Research and Development Program (2016YFB0900103) and the National Natural Science Foundation of China (51477082).
	H. Zhang, Z. Hu and Y. Song are with the Department of Electrical Engineering, Tsinghua University, Beijing, 100084, P.~R.~China.
	
	E. Munsing and S. J. Moura are with the Department of Civil and Environmental Engineering, University of California, Berkeley, California, 94720, USA.
	}
}

\maketitle

\begin{abstract}
	%Distributed Energy resources (DERs), e.g., battery storage systems, plug-in electric vehicles, may reap significant revenue by providing regulation services to the power systems. %Compared with traditional generators, their flexibilities are constrained by not only their power capacities but also energy capacities. %In other words, they have memories for their cumulative energy consumption. 
	This paper studies the behavior of a strategic aggregator offering regulation capacity on behalf of a group of distributed energy resources (DERs, e.g. plug-in electric vehicles) in a power market. Our objective is to maximize the aggregator's revenue while controlling the risk of penalties due to poor service delivery. To achieve this goal, we propose data-driven risk-averse strategies to effectively handle uncertainties in: 1) The DER parameters (e.g., load demands and flexibilities) and 2) sub-hourly regulation signals (to the accuracy of every few seconds). We design both the day-ahead and the hour-ahead strategies. 
	In the day-ahead model, we develop a two-stage stochastic program to roughly model the above uncertainties, which achieves computational efficiency by leveraging novel aggregate models of both DER parameters and sub-hourly regulation signals.
	In the hour-ahead model, we formulate a data-driven  distributionally robust chance-constrained program to explicitly model the aforementioned uncertainties. This program can effectively control the quality of regulation service based on the aggregator's risk aversion. Furthermore, it learns the distributions of the uncertain parameters from empirical data so that it outperforms existing techniques, (e.g. robust optimization or traditional chance-constrained programming) in both modelling accuracy and cost of robustness.
	Finally, we derive a conic safe approximation for it which can be efficiently solved by commercial solvers. Numerical experiments are conducted to validate the proposed method.
\end{abstract}


%\vspace{-2mm}

\section{{Introduction}}
\label{sec:introduction}
This paper proposes a risk-averse regulation capacity offering strategy for an aggregator of distributed energy resources (DERs) in power markets. This strategy handles uncertainties including DER parameters and regulation signals. It balances the trade-off between the aggregator's expected revenue and its risk of penalty for poor regulation performance.

With the rapid development of smart meters and advanced control technologies, DERs such as battery storage systems\cite{taylor2017price}, plug-in electric vehicles (PEVs)\cite{Neyestani}, and shapeable loads\cite{ZhaoyuWang2017}, promise the ability of providing flexibility services to power systems. These services include filling load valleys or arbitraging temporal differences in energy market prices\cite{pourbabak2016distributed,BoxunXu2017,meng2017coordinated}, facilitating integration of renewable energy sources \cite{Wang2015}, providing auxiliary services \cite{V2G_Regulation_Ko2016,Liuhui_EVregulation_2016,wenzel2017real}, etc.
Among these, regulation or load-following services are often the most lucrative in many power markets, e.g., PJM\cite{PJM_pilong2013pjm}. 
However, the energy and power capacities of DERs are comparatively small, while most regulation markets require high power capacities to access the market. Therefore, a large population of DERs usually have to jointly participate in a regulation market under the coordination of an aggregator. During market operations, the aggregator first evaluates its DERs' available regulation capacities, and then bundles the services from those DERs as a bid in the market. During real-time operations, the aggregator recieves regulation signals from the system operator, and is responsible for adjusting the power consumed by the DERs to follow these signals while respecting each DER's parameters, e.g., load demands and rated power and energy capacities etc.

The aggregator is faced with the difficult challenge of constructing a profit-maximizing bid, as accurately evaluating DERs' regulation capacities ahead-of-time is difficult, and market participants are penalized for inability to deliver on promised regulation capacity. This uncertainty arises from a number of sources: 
\begin{enumerate}
	\item \textit{Parameters of the DERs}. DERs' own parameters, e.g., rated energy and power capacities, customer demands of shapeable loads, plug-in time of PEVs etc., directly constrain their available regulation capacities but are usually stochastic.
	\item \textit{Regulation dispatch signals}. The limited energy capacities of DERs are easily saturated, e.g., PEVs may be fully charged or discharged after a period of regulation operations; thus biased regulation signals may significantly affect their regulation performance. Furthermore, frequently responding to regulation signals may lead to considerable energy losses due to inefficiency, which can also affect DERs' regulation performance.
\end{enumerate}

Both of these are difficult to forecast: DER parameters and regulation dispatch signals are ultimately shaped by outside factors such as weather and consumer behavior, making both of them fundamentally stochastic. Furthermore, regulation signals have high temporal granularity (e.g. every 2 seconds in the PJM market \cite{PJM_pilong2013pjm}) so that explicitly modeling them is computationally expensive.

Although utilizing DERs for regulation services has been an important research area for years, these difficulties have not been addressed in a way that effectively balances the revenue and risk for DERs providing  regulation services. The relevant methodologies proposed in published papers can be generally divided into four categories:
\subsubsection{Deterministic programming} These works assume that uncertain parameters can be accurately forecast, e.g., PEV driving patterns in \cite{V2G_Sortomme2012}, or adopt their expected values, e.g., regulation signals in \cite{V2G_Sortomme2012,Zhang_V2G_TPS2017,Zhang2017}. To avoid forecasting individual DER parameters, some papers aggregate parameters into virtual large-scale DERs which are more stable, e.g., PEVs in \cite{AggregateModelV2G_Xu2016,Zhang_V2G_TPS2017}, or thermostatic loads in \cite{AggregateModel_Mathieu2014}.
These approaches can only roughly estimate the regulation capacities, but they may be too optimistic for regulation capacity offering.
\subsubsection{Two-stage stochastic programming}
In this approach, a finite number of future scenarios are first generated based on forecasting or Monte Carlo simulation. Then, the aggregator uses scenario-based two-stage stochastic programming to estimate the regulation capacities in the future and create a market offer which has the best expected performance for all the given scenarios, e.g., \cite{PEV_Regulation_Vagropoulos2013,Storage_Regulation_He2016,Sarker2016}. 
However, to ensure adequate accuracy, the number of scenarios has to be large so that the problem may be computationally inefficient. Therefore, these works should adopt inaccurate approximation techniques, for example, references \cite{PEV_Regulation_Vagropoulos2013,Storage_Regulation_He2016,Sarker2016} only consider hourly average regulation signals. 

\subsubsection{Robust optimization} This approach pursues the optimal strategy when the ``worst-case scenario'' happens in the future. Yao et. al. \cite{V2G_Regulation_Yao2017} and Kazemi et. al. \cite{Kazemi2017}  apply robust optimization to handle hourly regulation signals, which are assumed to be bounded in predetermined intervals.  Vrettos et. al. \cite{Storage_Regulation_Vrettos2016} utilize robust optimization to describe sub-hourly regulation signals for commercial buildings.
The summation of the signals during an operation period are assumed to be bounded by a threshold.
Though this approach is usually computationally efficient, it may be unnecessarily conservative because the worst scenario rarely happens in practice.

\subsubsection{Risk-averse approach}
To overcome the limitations of the robust optimization approach, some papers propose to use risk-averse approaches, e.g., the chance-constrained programming or the conditional value of risk (CVaR). Their constraints are not required to be satisfied under the ``worst-case scenario'', but will be satisfied within certain (tunable) probability bounds. Vay\`{a} et. al. \cite{Vaya2015} adopt chance-constrained programming to model uncertain PEV driving behaviors. Yao et. al. \cite{V2G_Regulation_Yao2016} use the CVaR to describe the regulation revenue considering both uncertainties in PEV behaviors and in regulation prices. However, references \cite{Vaya2015,V2G_Regulation_Yao2016} do not model the risks associated with penalties from poor regulation service delivery. These papers also both use scenario-based approximations, introducing the limitations described above.

%Besides the aforementioned approaches, Cheng et. al. \cite{V2G_Regulation_Cheng2016} utilize Markov decision process to handle uncertainties of load demands, electricity prices, and sub-hourly regulation signals. They solve the problem by the dynamic programming which is computationally inefficient. This approach is a good technique for real-time operations but may be difficult to implement in ahead-of-time regulation capacity offering.

%\cite{Energycentric_Nosair2016} energy-based reserve definitions.

In this paper we advance this research by developing novel risk-averse data-driven regulation capacity offering strategies in both day-ahead and hour-ahead regulation markets for a DER aggregator. Compared with the aforementioned literature, the contributions of this paper are threefold:
\begin{enumerate}
	\item We formulate a two-stage stochastic programming model for day-ahead regulation capacity offering. This model adopts a novel hourly aggregate model to describe the regulation signals' influence on DERs' cumulative energy consumption, which is in small scale. Furthermore, it can accurately model DERs' charge and discharge inefficiency.
	\item We develop a risk-averse hour-ahead regulation capacity offering strategy based on the chance-constrained programming. The uncertainties of both the resource parameters and the sub-hourly regulation signals are explicitly modeled such that the trade-off between the revenue from providing regulation services and the risk of penalty for poor service delivery is effectively balanced. 
% 	\item We derive a continuous relaxation,  which is a standard linear chance-constrained program, for the aforementioned intractable hour-ahead model. We prove that the relaxation's solution is also feasible for the original problem with some mild assumptions. 
	\item We utilize historical market data to learn the information of the uncertain parameters' distributions. We then reformulate the hour-ahead program into a data-driven distributionally robust chance-constrained program based on the $\phi$-divergence. After that, we provide its convex robust counterpart in the form of second order cone programming (SOCP) so that it can be efficiently solved with off-the-shelf solvers.
	%\item We solve the aforementioned joint chance-constrained program by an efficient and scalable iterative algorithm based on the second order cone programming (SOCP) instead of scenario-based approximation.
\end{enumerate}
We validate the proposed strategy with numerical experiments using PEVs as an example. To the best of our knowledge, this is the first time that an hourly aggregate model for regulation signals has been designed and the first time that a data-driven distributionally robust chance-constrained programming is used for regulation capacity offering. 

The day-ahead two-stage stochastic strategies are introduced in Section II. Section III describes the risk-averse hour-ahead strategy, its relaxation, and its solution method. Numerical experiments are presented in Section IV. Section V concludes this paper.

\textbf{Notations:} We denote sets and functions by uppercase calligraphic English letters or uppercase Greek letters, e.g., $\mathcal{X}$, parameters by lowercase letters, e.g., $x$, and decision variables by uppercase English letters, e.g., $X$. We use boldface letters to represent vectors or matrices, e.g., $\mathbf{x}$ or $\mathbf{X}$, and mark stochastic parameters with the tilde sign, e.g., $\tilde{x}$.

\section{Day-ahead Regulation Capacity Offering}
This section describes a day-ahead regulation capacity offering strategy for an aggregator of  DERs, assuming the aggregator is a price-taker. The market environment is based on the PJM  market\cite{PJM_pilong2013pjm} but can be adapted to other markets. The nomenclature of this strategy is summarized in Table~\ref{tab_nomenclature1}.

The PJM energy market closes at 12:00 the day before the operating day, and the regulation market closes 60 minutes prior to the operating hour. Because the aggregator's power schedule in the energy market will affect its regulation capacities, we assume that the aggregator jointly schedules its day-ahead power profile and regulation capacities before the energy market closes. During the operating hours, the aggregator's actual power profile may deviate from its day-ahead offers, but the imbalances will be paid for based on real-time electricity prices. Its regulation capacity offers can also be reduced (but not increased) before the regulation market closes\cite{PJM_pilong2013pjm}. 

\begin{figure*}
	\centering
	\vspace{-8mm}
	\subfigure[Original trajectory]{
		\includegraphics[width=0.64\columnwidth]{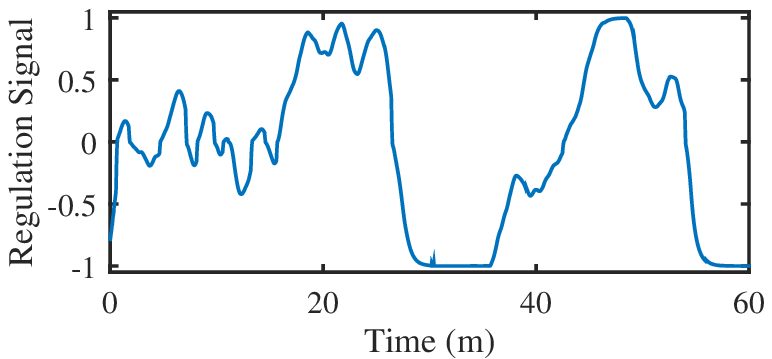}
		\label{fig_signal_original}
	}
	\subfigure[Arranged trajectory]{
		\includegraphics[width=0.64\columnwidth]{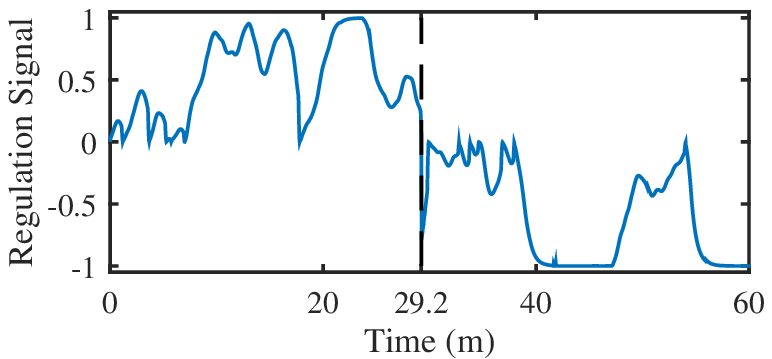}\label{fig_signal_arrange}
	}
	\subfigure[Hourly aggregate trajectory]{
		\includegraphics[width=0.64\columnwidth]{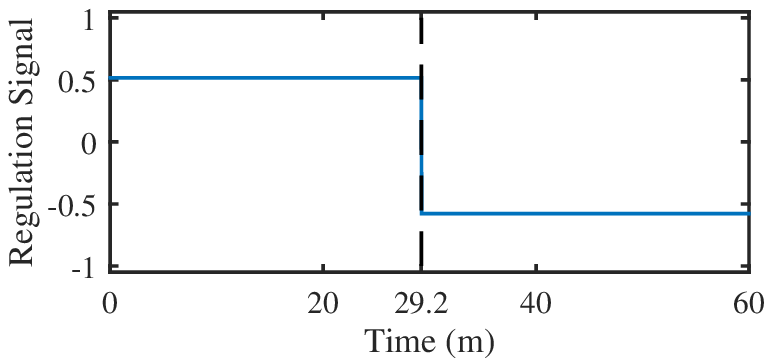}\label{fig_signal_average}
	}
	\vspace{-2mm}
	\caption{The trajectories of regulation signals, i.e., (a), $\mathbf{s}_t$, (b), $\mathbf{s}_t^{'}$, (c), $\mathbf{s}_t^{''}$.}
	\vspace{-5mm}
\end{figure*}
\subsection{Aggregate Model of Regulation Signals}\label{signal}
Frequently charging and discharging may lead to significant energy losses due to inefficiency. However, because the regulation signals are highly stochastic and their time granularity is very small, it is computationally expensive to explicitly model them in the capacity offering strategy.
This section proposes a novel hourly aggregate model to describe the regulation signals. Compared with published papers using hourly average values, this model can evaluate inefficiency with greater accuracy.

A sample regulation signal, $s_{t,d} \in [-1, 1]$, for one hour in the the PJM market is shown in Fig.~\ref{fig_signal_original}. The signals are issued every 2 seconds. When $s_{t,d}>0$, up regulation is required and the aggregator should decrease the charging power or increase the discharging power; when $s_{t,d}<0$ the opposite is true. The value of $s_{t,d}$ is the ratio of the aggregator's power shift (increment or decrement) to its offered capacity; we represent the trajectory as a vector $\mathbf{s}_t = \{s_{t,1}, s_{t,2}, ..., s_{t,1800}\}$. 

To develop the aggregate model for the regulation signals, we assume that the energy and power constraints are not binding during each sub-hourly interval. Then, the DERs' actual power during each sub-hourly interval $d$ is a function of their power baseline $P_{t}^\text{gr,da}$ and the signal ${s}_{t,d}$; we denote it by $\mathcal{P}(s_{t,d}, P_{t}^\text{gr,da})$. The cumulative energy consumption of the DERs is thus a function of $P_{t}^\text{gr,da}$ and the trajectory of ${s}_{t,d}$
and is denoted by $\mathcal{E}(\mathbf{s}_t, P_{t}^\text{gr,da})$. Then, we have $\mathcal{E}(\mathbf{s}_t, P_{t}^\text{gr,da}) = \sum_{d}\mathcal{P}(s_{t,d}, P_{t}^\text{gr,da})\Delta d$. 

If we think of the regulation signals not as an ordered vector but instead as a set of unordered data, we can use $\mathcal{S}_t$ to denote the set of regulation signal trajectories which can be created by rearranging the data in $\mathbf{s}_t$.
As an example, we may create $\mathbf{s}_t^{'} = \{s_{t,1}^{'}, s_{t,2}^{'}, ...,s_{t,d^{*}}^{'}, ..., s_{t,1800}^{'}\}\in \mathcal{S}_t$, in which, $ s_{t,d}^{'}\geq 0, \forall d\leq d^{*}, \text{and } s_{t,d}^{'}\leq 0, \forall d> d^{*}$. 
We can see that $\mathbf{s}_t^{'}$ can easily be created from $\mathbf{s}_t$ by moving all the regulation-up signals forward and the regulation-down signals backward, and thus confirm that $\mathbf{s}_t^{'}\in \mathcal{S}_t$. 
Fig.~\ref{fig_signal_arrange} shows an example of how this might be done with $d^{*}=726$. 
Hence, we have:
\begin{proposition}\label{pro_1}
	\textit{If the power and energy constraints are not binding, then }$\forall \mathbf{s}_t, \forall P_{t}^\text{gr,da}, \forall  \mathbf{s}_t^{'}\in \mathcal{S}_t$,  $\mathcal{E}(\mathbf{s}_t, P_{t}^\text{gr,da})$=$\mathcal{E}(\mathbf{s}_t^{'}, P_{t}^\text{gr,da})$. 
\end{proposition}
This proposition is intuitive, because $\mathcal{E}(\mathbf{s}_t, P_{t}^\text{gr,da}) = \sum_{d}\mathcal{P}(s_{t,d}, P_{t}^\text{gr,da})\Delta d= \sum_{d}\mathcal{P}(s_{t,d}^{'}, P_{t}^\text{gr,da})\Delta d=\mathcal{E}(\mathbf{s}_t^{'}, P_{t}^\text{gr,da})$.
In the day-ahead scheduling problem, we can assume that the aggregator can fulfill all the regulation requirements and relax the energy and power constraints so that Proposition \ref{pro_1} holds. 

Then, we can further approximate $\mathbf{s}_t^{'}$ by $\mathbf{s}_t^{''} = (s_{t,1}^{''}, s_{t,2}^{''}, ...,s_{t,d^{*}}^{''}, ..., s_{t,1800}^{''})$, in which, $ s_{t,d}^{''} =s_t^\text{up}= \frac{1}{d^{*}}\sum_{d=1}^{ d^{*}}{s_{t,d}^{'}}, \forall d\leq d^{*}, \text{and } s_{t,d}^{''} =s_t^\text{dn}= \frac{1}{1800-d^{*}}\sum_{d=d^{*}+1}^{ 1800}{s_{t,d}^{'}}, \forall d> d^{*}$ (see Fig.~\ref{fig_signal_average} for an example). Thus, we have $\mathcal{E}(\mathbf{s}_t, P_{t}^\text{gr,da}) \approx \mathcal{E}(\mathbf{s}_t^{''}, P_{t}^\text{gr,da})$ so that $\mathbf{s}_t^{''}$ can be used to describe the regulation signals' influence on cumulative power consumption and regulation capacities. 

$\mathbf{s}_t^{''}$ is the proposed aggregate model for the regulation signals. It has only four parameters including the average regulation up and down signals, i.e., $s_t^\text{up}$ and $s_t^\text{dn}$, and the corresponding up and down period, i.e., $\delta t_t^\text{up}=\frac{d^*}{1800}$ hour and $\delta t_t^\text{dn}=1- s_t^\text{up}$ hour. In Fig.~\ref{fig_signal_average}, $s_t^\text{up}=0.52$, $s_t^\text{dn}=-0.58$, $\delta t_t^\text{up}=29.2$ minutes and $\delta t_t^\text{dn}=30.8$ minutes. 
%\textcolor{red}{A limitation of the above aggregate model is that it can not reflect the information of regulation mileage, i.e., $m_t$. Therefore, we should also use historical data to forecast the future $m_t$.}

\begin{remark}
	Using the proposed hourly aggregate parameters of the signals, $s_t^\text{up/dn}$and $\delta t_t^\text{up/dn}$, the scale of the modeling is significantly decreased. Furthermore, we can also estimate the sub-hourly charging and discharging inefficiency because charged and discharged energy can be easily distinguished as is shown in the following section. 
\end{remark}

\subsection{Two-stage Stochastic Programming Strategy}
In the day-ahead problem, we propose using the two-stage stochastic programming to roughly incorporate uncertainties. A set of scenarios of the aggregate regulation signals, $s_{\omega,t}^\text{up/dn}$, the \textit{regulation mileages}\footnote{The regulation mileage is the absolute summation of movement requested by the regulation signal. It is used to evaluate the regulation service contribution of the regulation resources in the PJM market\cite{PJM_pilong2013pjm}.}, $m_{\omega,t}$, the energy and power capacities (at the resource side), $e_{\omega,t}^{+/-}$ and $p_{\omega,t}^{+/-}$, are generated based on historical data. The energy offer $P_t^\text{gr,da}$ is the first-stage decision variable, and the regulation capacity offer of each scenario $R_{\omega,t}$ is the second-stage decision variable.

\subsubsection{Day-ahead Objective}
The objective of the aggregator in the day-ahead market is:
\begin{align}
	\max ~ \mathcal{F}= & \sum_{\omega\in \Omega}\sum_{t}\pi_{\omega} \left(c_t^\text{rc}  + c_t^\text{rp} m_{\omega,t} \right) R_{\omega,t}\Delta t \notag\\
	&-\sum_{t}c_t^\text{e,da} P_t^\text{gr,da}\Delta t -\sum_{\omega\in \Omega}\sum_{t}\pi_{\omega}{c^\text{d}}E_{\omega,t}^{\text{der,d}} \label{eqn_obj1}
	%&-\sum_{t} \lambda \left(P_{t}^\text{gr,da}-\overline{P_{t}^\text{gr,da}}\right)^2.\label{eqn_obj1_4}
\end{align}
The first term in (\ref{eqn_obj1}) is the revenue for providing regulation services, which includes two parts: the revenue for the cleared regulation capacities, and the revenue for the service performance evaluated by the regulation mileage\cite{PJM_pilong2013pjm}. The second term is the energy purchase cost. For some DERs, e.g., batteries or PEVs, discharging may cause degradation. Therefore, we add the last term to calculate the degradation cost, which is assumed to be proportional to the aggregate discharged energy, $E_{\omega,t}^\text{der,d}$.
%The last term (\ref{eqn_obj1_4}) is adopted to penalized undesirable base load variations. Because a V2G fleet has fast response capability, its ramp rate maybe very high such that it may cause ramp problems for the grid. .

\subsubsection{Day-ahead Constraints}
The charging and discharging power must not violate the capacities of the DERs:
\begin{align}
&P_{t}^\text{gr,da}-s_{\omega,t}^\text{up/dn}R_{\omega,t} = \frac{P_{\omega,t}^\text{c,up/dn}}{\eta^\text{c}} + P_{\omega,t}^\text{d,up/dn}\eta^\text{d}, ~~ \forall \omega \in \Omega, \forall t,\label{eqn_1_1}\\
& 0 \leq P_{\omega,t}^\text{c,up/dn} \leq (1-D_{\omega,t}^\text{up/dn})p_{\omega,t}^{+}, \qquad \forall \omega \in \Omega, \forall t,\label{eqn_1_3}\\
& D_{\omega,t}^\text{up/dn}p_{\omega,t}^{-} \leq P_{\omega,t}^\text{d,up/dn} \leq 0, \qquad \forall \omega \in \Omega, \forall t,\label{eqn_1_4}\\
& D_{\omega,t}^\text{up}, D_{\omega,t}^\text{dn} \in \{0,1\}, \qquad \forall \omega \in \Omega, \forall t, \label{eqn_1_7}\\
&P_{\omega,t}^\text{der,up/dn}={P_{\omega,t}^\text{c,up/dn}} + P_{\omega,t}^\text{d,up/dn}, \qquad \forall \omega \in \Omega, \forall t,\label{eqn_1_5}\\
&E_{\omega,t}^\text{der}=\delta t_{\omega,t}^\text{up} P_{\omega,t}^\text{der,up} + \delta t_{\omega,t}^\text{dn} P_{\omega,t}^\text{der,dn}, \qquad \forall \omega \in \Omega, \forall t,\label{eqn_1_6}\\
&e_{\omega,t}^{-} \leq \sum_{\tau=t_0}^{t}{E_{\omega,\tau}^\text{der}} \leq e_{\omega,t}^{+}, \qquad \forall \omega \in \Omega, \forall t.\label{eqn_1_8}
\end{align}
Equation (\ref{eqn_1_1}) calculates the expected power during the sub-hourly interval when regulation up/down is required, i.e., $\delta t_{\omega,t}^\text{up/dn}$ in the aggregate regulation signal model. It reflects the relationship between the grid-side power and the resource-side power considering charging and discharging efficiency. Equations (\ref{eqn_1_3}) and (\ref{eqn_1_4}) constrain the charging and discharging power respectively; the two constraints will not be active simultaneously as can be seen from constraint (\ref{eqn_1_7}). Equation (\ref{eqn_1_5}) calculates the expected charging and discharging power at the resource side during each sub-hourly interval. Equation (\ref{eqn_1_6}) calculates the total energy consumption at the resource side in each hour. 
After each hour, the cumulative energy consumption of the DERs can not violate their aggregate upper and lower energy capacities, which are guaranteed by equation (\ref{eqn_1_8}).

The total discharged energy in each hour is:
\begin{align}
	&E_{\omega,t}^\text{der,d}=-\delta t_t^\text{up}P_{\omega,t}^\text{d,up} - \delta t_t^\text{dn}P_{\omega,t}^\text{d,dn}, \qquad \forall \omega \in \Omega, \forall t.\label{eqn_1_9}
\end{align}

 Because the aggregator may reduce its regulation capacity offer 60 minutes before the actual operating hour (without any penalty), it is acceptable to be aggressive in the regulation capacity offering. Therefore, the aggregator can submit the highest capacities to the market:
\begin{align}
&R_{t}^\text{da}=\max_{\omega \in \Omega}R_{\omega,t}, \qquad \forall t.\label{eqn_1_10}
\end{align}
\vspace{-2mm}
\begin{remark}\label{remark_2}The aforementioned model, i.e., (\ref{eqn_obj1})--(\ref{eqn_1_10}), is a mixed integer linear program, which can be efficiently solved by the Branch-and-Cut algorithm.
	Furthermore, the scale of the problem is moderate because it employs the aggregated hourly regulation signals developed in Section \ref{signal} and adopts the aggregated power and energy constraints of large-scale DERs whose scales are irrelevant to the DERs' population. (The introduction of the latter is omitted for brevity but can be found in published literature, e.g., \cite{AggregateModel_Mathieu2014,AggregateModelV2G_Xu2016,Zhang_V2G_TPS2017}) %As a result, the proposed model is computationally efficient.
\end{remark}

\begin{table}
	\renewcommand{\arraystretch}{1.05}
		\vspace{-2mm}
	\begin{footnotesize}
		\caption{Nomenclature of the day-ahead strategy }\label{tab_nomenclature1}
		\vspace{-1mm}
		\begin{tabular}{p{0.75cm}p{7.cm}}
			\hline
			\multicolumn{2}{c}{\textbf{Indices/sets}}\\
			$\omega/\Omega$&Index/set of all the scenarios.\\
			$d/\kappa$& Index of sub-hourly intervals, $\Delta$d=2 seconds.\\
			$t/\tau$& Index of hours, $\Delta$t=1 hour, $t_0$ is the initial hour.\\
			\hline
			\multicolumn{2}{c}{\textbf{Parameters}}\\
			$\eta^\text{c/d}$& Charge/discharge efficiency.\\
			$\pi_{\omega}$&Occurrence probability of scenario $\omega$.\\
			$c_t^\text{e,da}$&Day-ahead whole-sale energy price forecast, in \$/kWh.\\
			$c_t^\text{rc}$&Regulation capacity price forecast at time $t$, in \$/kWh.\\
			$c_t^\text{rp}$&Regulation performance price forecast at time $t$, in \$/kWh.\\
			$c^\text{d}$& Degradation costs due to discharging, in \$/kWh.\\
			$e_{\omega,t}^{+/-}$& Aggregate upper/lower energy capacity (at the resource side), in \$/kWh.\\
			$m_{\omega,t}$& Regulation millage. \\
			$p_{\omega,t}^{+/-}$& Aggregate upper/lower power capacity (at the resource side), in \$/kW.\\
			$\mathbf{s}_{t}$& Trajectory of regulation signals.\\
			$s_{t,d}$& Regulation signal.\\
			
			$s_t^\text{up/dn}$& Average regulation up/down signal in hour $t$.\\    
			$\delta t_t^\text{up/dn}$& Length of regulation up/down period in hour $t$, in hour.\\
			
			\hline
			\multicolumn{2}{c}{\textbf{Decision variables}}\\
			$D_{\omega,t}^\text{up/dn}$& Binary discharge decision variable: $D_{\omega,t}^\text{up/dn}=1$, if the DERs are discharged; $D_{\omega,t}^\text{up/dn}=0$, otherwise.\\
			$E_{\omega,t}^\text{der}$&The total energy consumption at the resource side, in kWh.\\
			$E_{\omega,t}^\text{der,d}$&The total discharged energy at the resource side, in kWh.\\
			
			$P_{\omega,t}^\text{der,up/dn}$&The power at the resource side when up/down regulation, i.e., $s_{\omega,t}^\text{up/dn}$ is dispatched, in kW.\\
			$P_{\omega,t}^\text{gr,da}$&The total power schedules at the grid side, in kW.\\
%			$P_{\omega,t}^\text{G,up/dn}$&The power at the grid side when up/down regulation, i.e., $s_{\omega,t}^\text{up/dn}$ is dispatched, in kW\\
			$P_{\omega,t}^\text{c/d,up/dn}$&The charging/discharging power at the resource side when up/down regulation, i.e., $s_{\omega,t}^\text{up/dn}$ is dispatched, in kW\\
			$R_{\omega,t}$ &Regulation capacity offer in each scenario, in kW\\
			$R_{t}^\text{da}$&Regulation capacity offer at day-ahead, in kW\\
			\hline
		\end{tabular}
	\end{footnotesize}
	\vspace{-4mm}
\end{table}

\section{Hour-ahead Regulation Capacity Offering}
In the hour-ahead problem, although the future available regulation capacities are uncertain, the aggregator still needs to give an explicit offer for the target hour.  If the aggregator fails to respond to the regulation signals as it offered, it may get penalized by the market. 
To control regulation service quality and balance revenue and risk, we propose a risk-averse strategy based on chance-constrained programming\cite{DCC_ben2009robust}. 

\subsection{Risk-averse Chance-constrained Strategy}
\subsubsection{Hour-ahead Objective}
The objective is to offer a proper regulation capacity to the market 60 minutes before the actual operation hour $t$ to maximize the overall expected revenue: 
\begin{align}
	\max \mathcal{F}= &(c_t^\text{rc} + c_t^\text{rp} \overline{m}_t) R_t\Delta t -c_t^\text{e,rt} \Delta P_t^\text{gr,ha}\Delta t \notag\\
	&-\sum_{\omega\in \Omega}\pi_{\omega}{c^\text{d}}E_{\omega,t}^{\text{der,d}}, \label{eqn_obj2}
\end{align}
This includes the regulation revenue (the first term), the energy purchase cost because of power deviations (the second term), and the cost for degradation (the last term).

To avoid myopic hour-ahead capacity offering, we adopt a receding horizon optimization to consider power schedules and regulation offers after $t$. The corresponding formulations are the same with those in the day-ahead strategy, and are omitted here for brevity. 

\subsubsection{Hour-ahead Constraints}
The actual regulation capacity offered to the market cannot exceed the day-head value:
\begin{align}
	&0 \leq  R_t \leq R_t^\text{da}\label{eqn_2_0}.
\end{align}

The absolute value of the power deviation is constrained by:
\begin{align}
	&\Delta P_t^\text{gr,ha}\geq P_t^\text{gr,ha}-P_t^\text{gr,da},\label{eqn_2_1}\\
	&\Delta P_t^\text{gr,ha}\geq -P_t^\text{gr,ha}+P_t^\text{gr,da}.\label{eqn_2_2}
\end{align}

The degradation costs are influenced by various factors, e.g., the power schedules, the regulation offers, and the signals. Because these costs are risk neutral, we still use the two-stage stochastic programming to estimate the expected value. Namely, $P_{\omega,t}^\text{der,d}$ should still satisfy the following constraints\footnote{The binary variables in (\ref{eqn_2_3}) are redundant and can be relaxed because simultaneously charging and discharging the DERs in this constraint can not increase the regulation capacity.}:
\begin{align}
&\text{(\ref{eqn_1_1})--(\ref{eqn_1_7}), (\ref{eqn_1_9})}.\label{eqn_2_3}
\end{align}

The DERs' charging and discharging power during each sub-hourly time interval should not violate their capacities: 
\begin{align}
&P_{t}^\text{gr,ha}-\tilde{s}_{t,d}R_{t} = {P_{t,d}^\text{c}}/{\eta^\text{c}} + P_{t,d}^\text{d}\eta^\text{d}, \qquad \forall d,\label{eqn_2_4}\\
& P_{t,d}^\text{c} \leq (1-D_{t,d})\tilde{p}_{t}^{+}, \qquad \forall d,\label{eqn_2_6}\\
& P_{t,d}^\text{d} \geq D_{t,d}\tilde{p}_{t}^{-} , \qquad \forall d,\label{eqn_2_7}\\
& D_{t,d} \in \{0,1\}, ~P_{t,d}^\text{c} \geq 0, ~P_{t,d}^\text{d} \leq 0, \qquad \forall d. \label{eqn_2_8}
\end{align}

Their cumulative energy consumption should also not violate its lower and upper capacities:
\begin{align}
%&P_{t,d}^\text{der}={P_{t,d}^\text{c}} + P_{t,d}^\text{d}, \qquad \forall d,\label{eqn_2_5}\\
&\tilde{e}_{0}  + \sum_{\kappa=1}^{d}{({P_{t,\kappa}^\text{c}} + P_{t,\kappa}^\text{d})}\Delta d \geq \tilde{e}_t^{-}, \qquad \forall d,\label{eqn_2_9}\\
&\tilde{e}_{0} + \sum_{\kappa=1}^{d}{({P_{t,\kappa}^\text{c}} + P_{t,\kappa}^\text{d})}\Delta d \leq \tilde{e}_t^{+}, \qquad \forall d,\label{eqn_2_10}
\end{align}
where, $\tilde{e}_{0}$ is the cumulative energy consumption after hour $t-1$ (before the operating hour $t$). Because the power schedule and regulation capacity offers at hour $t-1$ are already fixed, therefore, the aggregator can estimate the distribution of $\tilde{e}_{0}$ by simulations based on historical regulation signals. 

The above formulations form a mixed-integer linear regulation capacity offering strategy at hour-ahead:
\begin{align}
\textbf{P1:}\qquad &\max \quad \text{(\ref{eqn_obj2})} \qquad \text{s.t.:\quad (\ref{eqn_2_0})--(\ref{eqn_2_10}).}\notag
\end{align}

In \textbf{P1}, parameters $\tilde{\mathbf{\xi}} = \{\tilde{p}_{t}^{+/-}, \tilde{e}_{t}^{+/-}, \tilde{e}_{0}, \tilde{s}_{t,d}\}$ are all stochastic and may affect the aggregator's future regulation performance. Therefore, the aggregator should properly offer its regulation capacities according to the distribution of $\tilde{\mathbf{\xi}}$ so that it can reap adequate revenue and, at the same time, effectively control the quality of regulation services. 

Hence, we propose to use chance constraints to describe the DERs' energy and power limits, as follows:
\begin{align}
\textbf{P2:}\qquad &\max \quad \text{(\ref{eqn_obj2})} \qquad \text{s.t.:\quad (\ref{eqn_2_0})--(\ref{eqn_2_4}), (\ref{eqn_2_8}), and }\notag\\
&\text{Pr}_{\tilde{\mathbf{\xi}}\sim {\Xi}}\left\{ \text{(\ref{eqn_2_6})} \right\} \geq 1-\epsilon, \qquad \forall d, \label{eqn_cc1}\\
&\text{Pr}_{\tilde{\mathbf{\xi}}\sim {\Xi}}\left\{ \text{(\ref{eqn_2_7})} \right\} \geq 1-\epsilon, \qquad \forall d, \label{eqn_cc2}\\
&\text{Pr}_{\tilde{\mathbf{\xi}}\sim {\Xi}}\left\{ \text{(\ref{eqn_2_9})} \right\} \geq 1-\epsilon, \qquad \forall d, \label{eqn_cc3}\\
&\text{Pr}_{\tilde{\mathbf{\xi}}\sim {\Xi}}\left\{ \text{(\ref{eqn_2_10})} \right\} \geq 1-\epsilon, \qquad \forall d, \label{eqn_cc4}
\end{align}
which ensures that the probability that each of the power and energy constraint is violated (which means the aggregator fails to fully respond to the regulation signal) is less than $\epsilon$. The trade-off between the service revenue and quality can be effectively balanced by tuning $\epsilon$, which can be determined according to the aggregator's risk preference.

\begin{remark}\label{remark_3}
\textbf{P2} is intractable in its current form, because:
\begin{enumerate}
	\item The scale of binary variables is large due to the small time granularity, i.e, $\Delta d$. Furthermore, the binary variables can not be directly relaxed; otherwise, \textbf{P2} may overestimate the regulation capacities.\footnote{Without the binary constraints, the program may provide down regulations by charging and discharging the DERs simultaneously to waste electricity.}
	\item The probabilities in the left-hand sides of constraints (\ref{eqn_cc1})--(\ref{eqn_cc4}) are hard to evaluate. These are not convex even the binary variables are relaxed\cite{DCC_ben2009robust}.
	\item The decision variables $P_{t,d}^\text{c}$ and $P_{t,d}^\text{d}$ in the constraints (\ref{eqn_cc1})--(\ref{eqn_cc4}) are intermediate variables which are not explicit functions of the uncertain parameters $\tilde{\mathbf{\xi}}$. Hence, \textbf{P2} is not in a standard chance-constrained program form and cannot be directly solved by existing techniques.
\end{enumerate}
\end{remark}

In the following sub-sections, we first derive a safe continuous relaxation for \textbf{P2} and reformulate it into a typical chance-constrained program; then, we solve it as a \textit{data-driven distributionally robust chance-constrained program}\cite{DCC_Data_Driven_Jiang2016} based on its convex safe approximation in the form of an SOCP. 

\begin{table}
	\renewcommand{\arraystretch}{1.05}
		\vspace{-6mm}
	\begin{footnotesize}
		\begin{threeparttable}
		\caption{Nomenclature of the hour-ahead strategy }\label{tab_nomenclature2}
		\vspace{-1mm}
		\begin{tabular}{p{0.75cm}p{7.cm}}
			\hline
			\multicolumn{2}{c}{\textbf{Indices/sets}}\\
			$\tilde{\mathbf{\xi}}$& Set of uncertain parameters, $\tilde{\mathbf{\xi}} = \{\tilde{p}_{t}^{+/-}, \tilde{e}_{t}^{+/-}, \tilde{e}_{0}, \tilde{s}_{t,d}\}$.\\
			${\Xi}$& Distribution of $\tilde{\mathbf{\xi}}$.\\
			\hline
			\multicolumn{2}{c}{\textbf{Parameters}}\\
			$c_t^\text{e,rt}$&Real-time energy price forecast at time $t$, in \$/kWh.\\
			$\tilde{e}_{0}$ & Cumulative energy consumption after hour $t-1$.\\
			 $\epsilon$ &The tolerance parameter of the chance constraints.\\
			 $\overline{m}_t$ & The expected regulation mileage.\\
			\hline
			\multicolumn{2}{c}{\textbf{Decision variables}}\\
			$\Delta P_t^\text{gr,ha}$& Absolute values of the power deviations between day-ahead and hour-ahead schedules, in kW.\\
			$D_{t,d}$& Binary discharge decision variable.\\
			$P_{t}^\text{gr,ha}$ & Hour-ahead power schedule, in kW.\\
			$P_{t,d}^\text{c/d}$& Charging/discharging power at the resource side, in kW.\\
%			$P_{t,d}^\text{der}$& Total power at the resource side, in kW.\\
			$R_{t}$ &Actual (hour-ahead) regulation capacity offer, in kW.\\
			\hline
		\end{tabular}
		\begin{tablenotes}
			\item \textbf{Note:} We omitted those notations already appeared in Table~\ref{tab_nomenclature1}.
		\end{tablenotes}
	\end{threeparttable}
	\end{footnotesize}
	\vspace{-5mm}
\end{table}

\subsection{Relaxation of the Chance-constrained Program}
In this section, we relax the binary variables and eliminate the intermediate variables of \textbf{P2} to reformulate it into a typical linear chance-constrained program. Due to space limitation, we only provide the results. Interested readers can refer to the supplementary material for the detailed proof. 

Our derivations lead to the following proposition:
\begin{proposition}
	\textit{When $\tilde{p}_{t}^{+}\geq 0$, $\tilde{p}_{t}^{-}\leq 0$, and }$P_{t}^\text{gr,da}$ \textit{and} $P_{t}^\text{gr,ha}$ \textit{have the same sign, the constraints (\ref{eqn_2_4})--(\ref{eqn_2_10}) hold if the following constraints hold:}
	\begin{align}
	& \eta^\text{c} \left(P_{t}^\text{gr,ha}-\tilde{s}_{t,1}^\text{a}R_{t}\right) - \tilde{p}_{t}^{+} \leq 0, \label{eqn_cc1_3}\\
	& -\frac{1}{\eta^\text{d}} \left(P_{t}^\text{gr,ha}-\tilde{s}_{t,1}^\text{a}R_{t}\right) +\tilde{p}_{t}^{-} \leq 0,\label{eqn_cc2_3}\\
	&\left\{\begin{array}{l}
	- \eta^\text{c}d \Delta d P_{t}^\text{gr,ha} + \left( \frac{1+\eta^\text{c}\eta^\text{d}}{2\eta^\text{d}}\tilde{s}_{t,d}^\text{a} +\frac{1-\eta^\text{c}\eta^\text{d}}{2\eta^\text{d}} \right)d \Delta d R_{t}  \\
	-  \tilde{e}_{0} + \tilde{e}_t^{-} \leq 0,  \quad\forall d, \quad\textit{~if~} P_{t}^\text{gr,da}\geq 0,\\
	- \frac{1}{\eta^\text{d}}d \Delta d P_{t}^\text{gr,ha} + \left( \frac{1+\eta^\text{c}\eta^\text{d}}{2\eta^\text{d}}\tilde{s}_{t,d}^\text{a} +\frac{1-\eta^\text{c}\eta^\text{d}}{2\eta^\text{d}} \right)d \Delta d R_{t}  \\
	-  \tilde{e}_{0} + \tilde{e}_t^{-} \leq 0,  \quad\forall d, \quad\textit{~if~} P_{t}^\text{gr,da}< 0,
	\end{array}\right. \label{eqn_cc3_3}\\
	& \eta^\text{c}d \Delta d P_{t}^\text{gr,ha} - \eta^\text{c}\tilde{s}_{t,d}^\text{a} d \Delta d R_{t}  + \tilde{e}_{0} -\tilde{e}_t^{+} \leq 0, \quad\forall d, \label{eqn_cc4_3}
	\end{align}
	\textit{where, $\tilde{s}_{t,d}^\text{a}=\frac{1}{d}\sum_{\kappa=1}^{d}{\tilde{s}_{t,\kappa}}$ is the average regulation signals from sub-hourly interval 1 to $d$; $\tilde{s}_{t,1}^\text{a}= \tilde{s}_{t,1}$.}
\end{proposition}

We also observe that the energy constraints (\ref{eqn_cc3_3})--(\ref{eqn_cc4_3}) are more conservative when $d$ is large (see the supplementary material). Therefore, we need only retain the constraints when $d=1800$. As a result, we obtain a new problem:
\begin{align}
\textbf{P3:}\qquad &\max \quad \text{(\ref{eqn_obj2})} \qquad \text{s.t.:\quad (\ref{eqn_2_0})--(\ref{eqn_2_3}), and }\notag\\
&\text{Pr}_{\tilde{\mathbf{\xi}}^{'}\sim {\Xi}^{'}}\left\{ \text{(\ref{eqn_cc1_3})} \right\} \geq 1-\epsilon, \label{eqn_cc1_2}\\
&\text{Pr}_{\tilde{\mathbf{\xi}}^{'}\sim {\Xi}^{'}}\left\{ \text{(\ref{eqn_cc2_3})} \right\} \geq 1-\epsilon, \label{eqn_cc2_2}\\
&\text{Pr}_{\tilde{\mathbf{\xi}}^{'}\sim {\Xi}^{'}}\left\{ \text{(\ref{eqn_cc3_3})} \right\} \geq 1-\epsilon, \qquad d=1800, \label{eqn_cc3_2}\\
&\text{Pr}_{\tilde{\mathbf{\xi}}^{'}\sim {\Xi}^{'}}\left\{ \text{(\ref{eqn_cc4_3})} \right\} \geq 1-\epsilon, \qquad d=1800. \label{eqn_cc4_2}
\end{align}
The above chance-constraints are all linear with uncertain parameters $\tilde{\mathbf{\xi}}^{'} = \{\tilde{p}_{t}^{+/-}, \tilde{e}_{t}^{+/-}, \tilde{e}_{0}, \tilde{s}_{t,1}^\text{a}, \tilde{s}_{t,1800}^\text{a}\}\sim {\Xi}^{'}$.

\begin{remark}
	\textbf{P3} is a typical linear chance-constrained  program. It is a safe approximation of \textbf{P2} under the mild assumptions that $\tilde{p}_{t}^{+}\geq 0$, $\tilde{p}_{t}^{-}\leq 0$, and $P_{t}^\text{gr,da}$ has the same sign with $P_{t}^\text{gr,ha}$. In other words, when the aforementioned assumptions are true, \textbf{P3}'s solution is also feasible for \textbf{P2}.
\end{remark}

\subsection{Data-driven Distributionally Robust Chance-constrained Program and its SOCP Approximation}
Chance-constrained programming is generally intractable. An individual chance constraint can be equivalently reformulated into its convex counterpart in only few cases, e.g., when the uncertain parameters are Gaussian \cite{DCC_ben2009robust}.  Some researchers propose the distributionally robust chance-constrained programming which does not require that the uncertain parameters follow specific distributions. By contrast, they only utilize parts of the uncertain parameters' information, e.g., their supports\cite{DCC_Calafiore2006}, or moments \cite{DCC_Calafiore2006,DCC_J_Zymler2013}. As a compromise, their solutions should be satisfied for any potential distribution with the same known information. Though the results may be conservative, they are usually tractable, e.g., in the form of SOCP\cite{DCC_Calafiore2006} or semi-definite programming\cite{DCC_J_Zymler2013}. 

Considering that historical regulation data in many power markets are public, we propose a data-driven approach to utilize the available information for the uncertain regulation signals. In this paper, we take the Reg-A type regulation signals in the PJM market as an example. Based on our analysis (given in the supplementary material), we observe that the signals in every two seconds, i.e., $\tilde{s}_{t,1}^\text{a}$ in (\ref{eqn_cc1_2})--(\ref{eqn_cc2_2}), do not follow a specific tractable distribution. By contrast, the hourly average regulation signals, i.e., $\tilde{s}_{t,1800}^\text{a}$ in (\ref{eqn_cc3_2})--(\ref{eqn_cc4_2}), are approximately Gaussian (hence, the hourly cumulative energy consumption $\tilde{e}_{0}$ is also approximately Gaussian).

Utilizing the above analysis and assuming that the forecast errors of the energy and power capacities are Gaussian, we adopt the distributionally robust chance-constrained programming (with known mean and covariance) to approximate the power constraints (\ref{eqn_cc1_2})--(\ref{eqn_cc2_2}). This approach only uses the empirical mean and covariance of the uncertain parameters. However, we utilize the Gaussian distribution to approximate the energy constraints (\ref{eqn_cc3_2})--(\ref{eqn_cc4_2}) (whose uncertain parameters are all approximately Gaussian). As the hourly regulation signals are not exactly Gaussian, we utilize the historical data to learn the $\phi$-divergence\cite{DCC_Data_Driven_Jiang2016} (a measure of the difference between two distributions) between the empirical distribution with the approximated Gaussian distribution. Then, we adopt the $\phi$-divergence-based data-driven distributionally robust chance-constrained programming\cite{DCC_Data_Driven_Jiang2016} to conservatively model the energy constraints. 

Before giving the further formulations, we first briefly introduce the $\phi$-divergence, a function that measures the distance between two nonnegative vectors $\mathbf{p}=(p_1, p_2, ..., p_n)^\intercal$ and $\mathbf{q}=(q_1, q_2, ..., q_n)^\intercal$. In this paper, we use it to measure the distance between two discrete probability distributions. We let $\mathbf{q}$ denote the true distribution's probabilities, and $\mathbf{p}$ denote the corresponding observations, so that we have $\sum_{i=1}^{n}{p_i}=\sum_{i=1}^{n}{q_i}=1$. The $\phi$-divergence between $\mathbf{p}$ and $\mathbf{q}$ is defined as follows: 
\begin{align}
&I_{\phi}(p, q)=\sum_{i=1}^{m}q_i\phi\left({p_i}/{q_i}\right),\label{eqn_phi}
\end{align}
in which, $\phi(t)$ is called the $\phi$-divergence function which is convex for $t\geq 0$, $\phi(1)=0$, $\phi(a/0)\doteq a\lim_{t\rightarrow\infty}\phi(t)/t$ for $a>0$, and $\phi(0/0)=0$. There are different types of $\phi$-divergence functions studied in published literature (see \cite{DCC_Data_Driven_Jiang2016}).  We adopt the $\chi^2$-divergence function, i.e., $\phi(t)=(t-1)^2$. 
We also use $I_{\phi}(\mathcal{D}, \mathcal{D}_0)$ to denote the $\phi$-divergence, i.e., $\chi^2$-divergence, between two distributions $\mathcal{D}$ and $\mathcal{D}_0$ (we can discretize the continuous distributions and calculate the $\phi$-divergence according to (\ref{eqn_phi})).

Letting $\mathbf{X}\doteq[R_t,  P_t^\text{gr,ha} ]^{\intercal}$ denote the vector of decision variables; $\hat{\mathbf{X}}\doteq[\mathbf{X}^{\intercal}~ 1]^{\intercal}$, \textbf{P3} can be reformulated into a data-driven distributionally robust chance-constrained program:
\begin{align}
&\textbf{P4:}~\max \quad \text{(\ref{eqn_obj2})} \qquad \text{s.t.:\quad (\ref{eqn_2_0})--(\ref{eqn_2_3}), and }\notag\\
&\inf_{\mathcal{D}_j \in (\overline{\mathbf{d}}_j, \mathbf{\Gamma}_j)}\text{Pr}_{\tilde{\mathbf{d}}_j\sim \mathcal{D}_j}\left\{ \tilde{\mathbf{d}}_j^{\intercal}\hat{\mathbf{X}}\leq 0 \right\} \geq 1-\epsilon,  ~~\forall j =1, 2, \label{eqn_CC3}\\
&\inf_{ I_{\phi}(\mathcal{D}_j, \mathcal{D}_0)\leq \rho}\text{Pr}_{\tilde{\mathbf{d}}_j\sim \mathcal{D}_j}\left\{ \tilde{\mathbf{d}}_j^{\intercal}\hat{\mathbf{X}}  \leq 0 \right\} \geq 1-\epsilon, ~\forall  j=3,4. \label{eqn_CC4}
\end{align}
Where $j$ is the index of chance-constraints, (\ref{eqn_CC3}) are the power constraints and (\ref{eqn_CC4}) are the energy constraints, $\tilde{\mathbf{d}}_j \doteq [\tilde{\mathbf{a}}_j^{\intercal},  \tilde{b}_j]^{\intercal}$ with mean $\overline{\mathbf{d}}_j$ and covariance $\mathbf{\Gamma}_j \doteq \text{var}\{{\tilde{\mathbf{d}}_j}\}\succeq 0$, and the entries of $\tilde{\mathbf{a}}_j$ and $\tilde{b}_j$ are linear combinations of the uncertain parameters in $\tilde{\mathbf{\xi}}^{'}$. We use $(\overline{\mathbf{d}}_j, \mathbf{\Gamma}_j)$ to denote the set of distributions having the same mean $\overline{\mathbf{d}}_j$ and covariance $\mathbf{\Gamma}_j$. 
$I_{\phi}(\mathcal{D}_j, \mathcal{D}_0)$ is the $\chi^2$-divergence between distribution $\mathcal{D}_j$ and $\mathcal{D}_0=Gaussian(\overline{\mathbf{d}}_j, \mathbf{\Gamma}_j)$, $\forall j=3, 4$; $\rho$ is the empirical upper bound of the $\chi^2$--divergence learned from historical data. 
 $\overline{\mathbf{d}}_j$ and $\mathbf{\Gamma}_j$ can also be easily calculated based on the empirical means and covariances of the parameters of $\tilde{\mathbf{\xi}}^{'}$. 

\begin{remark}
	\textbf{P4} is a safe approximation of the original chance-constrained program \textbf{P3}, i.e., the solution of \textbf{P4} is also feasible for \textbf{P3}. This is because the true distributions of the uncertain parameters in \textbf{P3} are included in the set $(\overline{\mathbf{d}}_j, \mathbf{\Gamma}_j)$, $\forall j=1, 2$, and set $\left\{\mathcal{D}_j:  I_{\phi}(\mathcal{D}_j, \mathcal{D}_0)\leq \rho \right\}$, $\forall j=3, 4$. In words, the power constraints (\ref{eqn_cc1_2})--(\ref{eqn_cc2_2}) are satisfied for all the distributions with the same mean $\overline{\mathbf{d}}_j$ and covariance $\mathbf{\Gamma}_j$ while the energy constraints (\ref{eqn_cc3_2})--(\ref{eqn_cc4_2}) are satisfied for all the distributions whose $\chi^2$--divergence between $\mathcal{D}_0=Gaussian(\overline{\mathbf{d}}_j, \mathbf{\Gamma}_j)$ is less than $\rho$.
\end{remark}

Based on the findings in \cite{DCC_Calafiore2006} and \cite{DCC_Data_Driven_Jiang2016}, we have:

\begin{proposition} \textit{For any $\epsilon \in (0,0.5]$, the chance constraints (\ref{eqn_CC3})--(\ref{eqn_CC4}) hold if the following constraints hold:} 
\begin{align}
&\sqrt{{(1-\epsilon)}/{\epsilon}} \| \mathbf{\Gamma}_j^{\frac{1}{2}} \hat{\mathbf{X}}\|_2  + \overline{\mathbf{d}}_j^{\intercal} \hat{\mathbf{X}} \leq 0, \qquad \forall j=1,2,\label{eqn_CC3_2}\\
&\Psi_G^{-1}(1-\epsilon^{'}
) \| \mathbf{\Gamma}_j^{\frac{1}{2}} \hat{\mathbf{X}}\|_2  + \overline{\mathbf{d}}_j^{\intercal} \hat{\mathbf{X}} \leq 0, \qquad \forall j=3,4,\label{eqn_CC4_2}
\end{align}
\textit{where, $\Psi_G(\epsilon)=(1/\sqrt{2\pi})\int_{-\infty}^{\epsilon}\text{exp}(\frac{t^2}{2})\text{d}t$; $\epsilon^{'}=\epsilon-\frac{\sqrt{\rho^2 +4\rho(\epsilon^2-\epsilon)}-(1-2\epsilon)\rho}{2\rho+2}$ is an adjusted tolerance based on $\chi^2$-divergence. }
\end{proposition}
Based on Proposition 3, we obtain the safe convex approximation of \textbf{P4}, which is an SOCP and can be efficiently solved:%\footnote{With the receding horizon formulations considering future hours' regulation capacity offering, \textbf{P5} is a mixed-integer SOCP. It can still be efficiently solved because the scale of the binary variables is small (see Remark \ref{remark_2}).}
\begin{align}
\textbf{P5:}\qquad &\max \quad \text{(\ref{eqn_obj2})} \qquad \text{s.t.:\quad (\ref{eqn_2_0})--(\ref{eqn_2_3}), and (\ref{eqn_CC3_2})--(\ref{eqn_CC4_2})}.\notag
\end{align}
%This program can be efficiently solved by an off-the-shelf commercial solver.

For brevity we only provide derivations for the formulation of constraints (\ref{eqn_CC3})--(\ref{eqn_CC4}), the calculation of $\overline{\mathbf{d}}_j$ and $\mathbf{\Gamma}_j$ (based on the parameters of $\tilde{\xi}^{'}$), and the proof for the Proposition 3 in the supplementary material.

%\vspace{-2mm}
\section{Experiments}
%\vspace{-2mm}
%\subsection{Parameter Settings}
We model a fleet of 5,000 private PEVs to verify the effectiveness of the proposed method. Only home charging/discharging is considered in this case. The Nissan Leaf PEV is chosen to represent the PEV population, with a battery capacity 24 kWh. We assume that 50\% of the customers install Level 1 chargers with +/-3.3 kW rated charging/discharging power. The others install Level 2 chargers with +/-6.6 kW rated power. The efficiencies ($\eta^{\text{c/d}}$) are 92\%; The battery degradation cost is assumed to be \$4.1 per charge cycle, i.e., $c^\text{b}$=4.1/24 \$/kWh. We adopt the method proposed in \cite{Zhang_V2G_TPS2017} to generate the PEVs' driving behaviors and the forecasted aggregate power and energy capacities. The standard variances of both the energy and power capacities are assumed to be 0.05. %The other PEV parameters were set the same with \cite{Zhang_V2G_TPS2017}.

We use average values of 28-day historical data as the forecasted prices of the energy and regulation market. The RegA-type regulation signals from Aug. 1st, 2015, to Jul. 31th, 2016 of the PJM market are used in the simulation to estimate their mean, covariance, $\chi^2$-divergence parameter $\rho$, and hourly mileages\cite{PJM_pilong2013pjm}. We conduct simulations, including regulation capacity offering and real-time operations, for 31 days to validate the performance of the proposed strategies, in which the real regulation signals in August 2015 are used as the input and the energy and power boundaries are generated by the Monte-carlo simulation.
We use CPLEX to solve the problems on a laptop with a 4-core Intel Core i7 processor and 8 GB memory. It takes a few seconds to solve each problem.

%\vspace{-2mm}
\subsection{Trade-off Between Revenue and Risk}
We first conduct experiments to validate the effectiveness of the proposed strategy at balancing the trade-off between the regulation revenue and risk of poor service delivery. The PJM market compensates the regulation DERs based on the cleared regulation capacities, the mileages and the DERs' performance scores. The last factor is used to reflect the accuracy of the DERs following the regulation signals. It consists of three parts: precision, correlation, and delay\cite{PJM_pilong2013pjm}. Considering that the simulations cannot truly reflect the real-time operation results, e.g., the delays, we use the precision score to approximate the performance score, which is calculated as follows:
\begin{align}
&S=100\%-\frac{1}{n}\sum_{d=1}^{n}\left|{\left(s_d - s_d^\text{r}\right)}/{\overline{|s|}}\right|,
\end{align}
where, $s_d$ is the instructed signal, $s_d^\text{r}$ is the actual response, $\overline{|s|}$ is the average of the absolute values of $s_d$, $n$ is the number of samples. We therefore can calculate the actual total revenue, i.e., $Revenue_\text{A}$, by adjusting its expected value (the objective of \textbf{P5}) with the performance score, as follows:
\begin{align}
Revenue_\text{A} = S\times Revenue_\text{R}-Cost_\text{der} - Cost_\text{D},
\end{align}
in which, $Revenue_\text{R}$ is the expected regulation revenue, $Costs_\text{der}$ is the energy cost, $Costs_\text{D}$ is the degradation cost. The expected total revenue is equal to $Revenue_\text{A}$ if $S=1$.

\begin{figure}
	\centering
	\vspace{-4mm}
	\begin{minipage}[t]{1\columnwidth}
		\centering
		\epsfig{file = 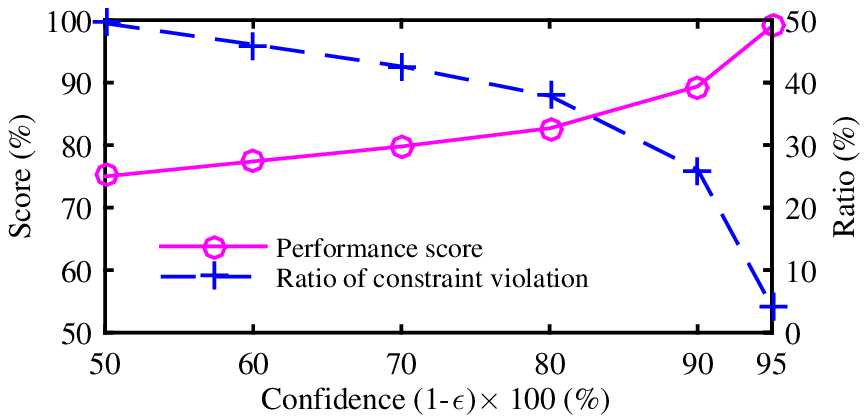, width = 7cm}
		\vspace{-3mm}
		\caption{Performance score and ratio of constraint violations.}
		\label{fig_score}
		\vspace{1mm}
	\end{minipage}
	\begin{minipage}[t]{1\columnwidth}
		\centering
		\epsfig{file = 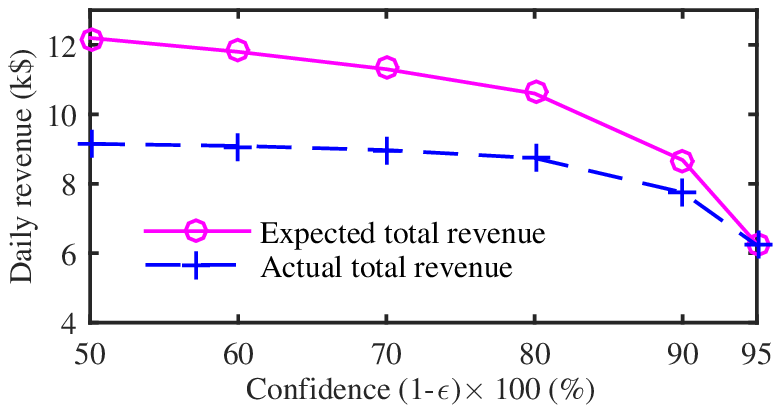, width = 7cm}
		\vspace{-3mm}
		\caption{Expected and actual daily total revenue.}
		\label{fig_benefit}
	\end{minipage}
	\vspace{-4mm}
\end{figure}

We conduct experiments under different tolerance levels, i.e., taking values of $1-\epsilon$ from 50\% to 95\%. The performance scores and the corresponding ratios of chance-constraint violations (unfulfilled regulation services) in the experiments are illustrated in Fig.~\ref{fig_score}. The average expected daily total revenue and the corresponding actual values are plotted in Fig.~\ref{fig_benefit}.

When $\epsilon$ increases, the probability that the regulation signals can not be precisely followed increases while the performance score decreases. With larger $\epsilon$, the aggregator will be more aggressive offering more regulation capacity to increase expected daily total revenue. However, the actual revenue growth is offset by the performance score when $\epsilon$ is larger than 30\%. This phenomena demonstrates the trade-off between revenue and the service quality (risk of poor service delivery). 

Note that the trade-off effect is highly dependent on the penalty mechanism. If we adopt the calculation method used in \cite{V2G_Regulation_Yao2017}, the performance score would be approximately equal to 100\% minus the ratio of chance-constraint violations, which will be much smaller (see Fig.~\ref{fig_score}). As a result, the trade-off would be more apparent. A review of different penalty mechanisms can be found in \cite{V2G_Regulation_Ko2016}.

%\vspace{-2mm}
\subsection{Benchmark with Other Strategies}

We adopt three benchmark strategies to validate the advantages of the proposed method: 
a) A robust strategy (Robust), which ensures that the regulation services can be delivered under the worst scenario. 
b) A deterministic strategy (Determ) which adopts expected values of the uncertain parameters. It is also used as a benchmark in \cite{V2G_Regulation_Yao2017}. 
%c) A deterministic strategy (Deter\_2) which assumes that the regulation signals have no effect on the DERs' cumulative energy consumption and requires that the regulation offers can be committed under the extreme signals, i.e., $\pm 1$. This is the strategy used in \cite{Zhang_Plan_TPS2016}. 
c) A risk-averse strategy (IgnoreEffi) which adopts the proposed method but ignores charging and discharging inefficiency at the capacity offering stage. After the capacities are offered, we simulate the real-time operations considering inefficiency to validate its performance. In the two risk-averse strategies (Proposed and IgnoreEffi), we both set $\epsilon=20\%$. 
The summary of their average daily regulation offers, performance scores, and actual total revenue are listed in Table~\ref{tab:case}. 

From the results, we can conclude that the robust approach is unnecessarily conservative with total daily regulation capacity offer only about 6\% of that in the proposed strategy- this clearly does not take full advantage of the DER flexibility. 

When assuming the uncertain parameters take average values, the aggregator may be very aggressive in offering regulation capacity because it is risk neutral. In fact, it provides regulation offers that are higher than its actual capacities (similar results were also found in \cite{V2G_Regulation_Ko2016}). However, the performance scores are significantly decreased with the average value falls below 70\%. As a result, its total benefit is lower than that of the proposed strategy. Furthermore, in the PJM market, if the performance scores are frequently below 75\%, it may fail the market's qualification test and become not eligible to participate in the market\cite{PJM_pilong2013pjm}.

When inefficiency is ignored, the aggregator may evaluate its regulation capacities less accurately. In Table~\ref{tab:case}, though all the experimental parameters (except the efficiency) are the same, the strategy IgnoreEffi's regulation capacity offer, performance score, and total revenue are all less than the proposed strategy. Compared with the proposed strategy, the total revenue is reduced by 3.8\%. An interesting point is that though IgnoreEffi offers less regulation capacity, its performance score is still lower that the proposed strategy. That is because the inefficiencies cause the DERs to consume additional electricity during real-time operations, affecting their regulation performance.

\begin{table}
	\renewcommand{\arraystretch}{1.}
	%\linespread{0.9}
	\vspace{-4mm}
	\centering
	\begin{small}
		\caption{Summary of different strategies}
		\vspace{-2mm}
		\begin{tabular}{cccc}
			\hline
			\multirow{2}{*}{Strategy}&{Regulation offer} & Performance &Total revenue \\
			&(MWh/day) & score (\%)& (k\$/day) \\
			\hline
			Proposed& 484 & 83 & 8.74\\
			Robust& 20.5 & 100 & 0.22 \\
			Determ& 511 & 68 & 7.96\\
			IgnoreEffi& 481 & 80 & 8.41 \\
			\hline
		\end{tabular}\label{tab:case}
	\end{small}
	\vspace{-4mm}
\end{table}

\section{{Conclusions}}
\label{sec:conclusion}
DERs are promising flexibility service providers for future power systems. We study strategic regulation capacity offering of an aggregator on behalf of a group of DERs in a power market. 
We propose a two-stage stochastic program for the day-ahead market. In this strategy, we design a novel aggregate model for sub-hourly regulation signals. Together with the aggregate model for DER parameters in published literature, it enhances both the computational efficiency and modeling accuracy of the proposed day-ahead strategy.
We formulate a data-driven distributionally robust chance-constrained program for the hour-ahead market which models uncertain parameters in greater granularity. This strategy can effectively balance the aggregator's regulation revenue and risk of poor regulation service delivery by tuning a simple parameter according to the aggregator's risk preference. We further derive an SOCP approximation for it so that it can be efficiently solved by a commercial solver. The proposed strategy outperforms published ones, e.g., those based on robust optimization or scenario-based risk-averse approaches, in both decision conservativeness and computational efficiency. 

We plan on studying the impact of of price uncertainty on optimal regulation offering strategies in future work.
%
%\appendices
%\section{}\label{appendix_1}

\bibliographystyle{ieeetr}
\bibliography{reference}

\end{document}

\title{
	Supplementary Material for ``Data-driven Chance-constrained Regulation Capacity Offering for Distributed Energy Resources"
	%Under Dispatch Uncertainty
	}

\author{
Hongcai~Zhang,~\IEEEmembership{Student Member,~IEEE,}
Zechun~Hu,~\IEEEmembership{Senior Member,~IEEE,}
Eric~Munsing,~\IEEEmembership{Student Member,~IEEE,}
Scott~J.~Moura,~\IEEEmembership{Member,~IEEE,}
and~Yonghua~Song,~\IEEEmembership{Fellow,~IEEE}
}

\maketitle

%\begin{IEEEkeywords}
%Energy-constrained resource, regulation service, risk-averse strategy, distributinally robust chance-constraint.
%\end{IEEEkeywords}

%\vspace{-2mm}

\section{Proof of Proposition 2}
Proposition 2 provides the safe continuous approximation for the following constraints:
\begin{align}
	&P_{t}^\text{gr,ha}-\tilde{s}_{t,d}R_{t} = \frac{P_{t,d}^\text{c}}{\eta^\text{c}} + P_{t,d}^\text{d}\eta^\text{d}, \qquad \forall d,\label{eqn_2_4}\\
	& P_{t,d}^\text{c} \leq (1-D_{t,d})\tilde{p}_{t}^{+}, \qquad \forall d,\label{eqn_2_6}\\
	& P_{t,d}^\text{d} \geq D_{t,d}\tilde{p}_{t}^{-} , \qquad \forall d,\label{eqn_2_7}\\
	& D_{t,d} \in \{0,1\}, ~P_{t,d}^\text{c} \geq 0, ~P_{t,d}^\text{d} \leq 0, \qquad \forall d, \label{eqn_2_8}\\
	&\tilde{e}_{0}  + \sum_{\kappa=1}^{d}{({P_{t,\kappa}^\text{c}} + P_{t,\kappa}^\text{d})}\Delta d \geq \tilde{e}_t^{-}, \qquad \forall d,\label{eqn_2_9}\\
	&\tilde{e}_{0} + \sum_{\kappa=1}^{d}{({P_{t,\kappa}^\text{c}} + P_{t,\kappa}^\text{d})}\Delta d \leq \tilde{e}_t^{+}, \qquad \forall d,\label{eqn_2_10}
\end{align}

\subsubsection{Power Capacity Constraints}We first derive the relaxation of the power capacity constraints (\ref{eqn_2_6})--(\ref{eqn_2_8}).
 Because when $\tilde{p}_{t}^{+}$ and $\tilde{p}_{t}^{-}$ have the same sign, e.g., both positive for shapeable load, the sign of the power in (\ref{eqn_2_4}) can be determined a priori. As a result, we do not need to distinguish charging and discharging power in (\ref{eqn_2_6})--(\ref{eqn_2_7}) so that they are reduced to continuous constraints. Here, we consider the scenario when $\tilde{p}_{t}^{+}\geq 0$ and $\tilde{p}_{t}^{-}\leq 0$, which is true for batteries or PEVs with vehicle-to-grid technology. In that case, constraint (\ref{eqn_2_6}) will be binding only when charging power happens and  $P_{t}^\text{gr,ha}-\tilde{s}_{t,d}R_{t}= {P_{t,d}^\text{c}}/{\eta^\text{c}}$. Therefore, it is equivalent to:
\begin{align}
	& \eta^\text{c} \left(P_{t}^\text{gr,ha}-\tilde{s}_{t,d}R_{t}\right) \leq \tilde{p}_{t}^{+}, \qquad \forall d, \label{eqn_4_2}
\end{align}
Similarly, constraint (\ref{eqn_2_7}) is binding only when discharging power happens so that it is equivalent to:
\begin{align}
	& \frac{1}{\eta^\text{d}} \left(P_{t}^\text{gr,ha}-\tilde{s}_{t,d}R_{t}\right) \geq \tilde{p}_{t}^{-}, \qquad \forall d.\label{eqn_4_3}
\end{align}
The above two constraints are both continuous and linear. Because the distribution of regulation signals in different time interval $d$ are the same, therefore, constraints (\ref{eqn_4_2})--(\ref{eqn_4_3}) for different $d$ are equivalent.

\subsubsection{Energy Capacity Constraints} Then, we then relax the energy capacity constraints (\ref{eqn_2_9})--(\ref{eqn_2_10}).
Using equation (\ref{eqn_2_4}) to eliminate the intermediate variable $P_{t,d}^\text{c}$, we can reformulate the  first energy constraint (\ref{eqn_2_9}) as follows:
\begin{align}
	&\tilde{e}_{0}+ \sum_{\kappa=1}^{d} \left(\eta^\text{c}\left(P_{t}^\text{gr,ha}-\tilde{s}_{t,\kappa}R_{t}\right) + (1-\eta^\text{c}\eta^\text{d}) P_{t,\kappa}^\text{d}\right)\Delta d\notag\\
	&\geq  \tilde{e}_t^{-},\qquad \forall d,\label{eqn_4_6}\\
	&\text{where,}~P_{t,\kappa}^\text{d} = \min\left(\frac{P_{t}^\text{gr,ha} -\tilde{s}_{t,\kappa}R_{t}}{\eta^\text{d}}, 0\right)\notag\\
	&\qquad\qquad~~= \frac{P_{t}^\text{gr,ha} -\tilde{s}_{t,\kappa}R_{t} - |P_{t}^\text{gr,ha}-\tilde{s}_{t,\kappa}R_{t}|}{2\eta^\text{d}} \notag\\
	&\qquad\qquad~~\geq \frac{P_{t}^\text{gr,ha}-|P_{t}^\text{gr,ha}| -(\tilde{s}_{t,\kappa}+1)R_{t}  }{2\eta^\text{d}}, \quad\forall \kappa. \label{eqn_4_7}
\end{align}

Because the hour-ahead power schedule, i.e., $P_{t}^\text{gr,ha}$, is close to the day-ahead power schedule, i.e., $P_{t}^\text{gr,da}$, we can safely assume that they have the same sign. As a result, the absolute value sign in (\ref{eqn_4_7}) can be eliminated a priori according to $P_{t}^\text{gr,da}$. Hence, if $P_{t}^\text{gr,da}\geq 0$, we can further eliminate variable $P_{t,d}^\text{d}$ based on (\ref{eqn_4_7}) to approximate constraints (\ref{eqn_4_6})--(\ref{eqn_4_7}) as:
\begin{align}
	&\tilde{e}_{0}+ \sum_{\kappa=1}^{d} \bigg(\eta^\text{c}\left(P_{t}^\text{gr,ha}-\tilde{s}_{t,\kappa}R_{t}\right) -\notag\\
	&\frac{1-\eta^\text{c}\eta^\text{d}}{2\eta^\text{d}} \left(\tilde{s}_{t,\kappa}R_{t}+R_{t}\right)\bigg)\Delta d \geq \tilde{e}_t^{-}, \qquad \forall d.\label{eqn_4_8_1}
\end{align}
Otherwise, if  $P_{t}^\text{gr,da}< 0$, we can approximate constraints (\ref{eqn_4_6})--(\ref{eqn_4_7}) as:
\begin{align}
&\tilde{e}_{0}+ \sum_{\kappa=1}^{d} \bigg(\eta^\text{c}\left(P_{t}^\text{gr,ha}-\tilde{s}_{t,\kappa}R_{t}\right) +\notag\\
&\frac{1-\eta^\text{c}\eta^\text{d}}{2\eta^\text{d}} \left(2P_{t}^\text{gr,ha} - \tilde{s}_{t,\kappa}R_{t}-R_{t}\right)\bigg)\Delta d \geq \tilde{e}_t^{-}, \quad \forall d.\label{eqn_4_8_2}
\end{align}

Because, $0\leq \eta^\text{c}, \eta^\text{d} \leq 1, P_{t,\kappa}^\text{d}\leq 0$, we also have
\begin{align}
	{P_{t,\kappa}^\text{c}} + P_{t,\kappa}^\text{d}&=\eta^\text{c}\left(P_{t}^\text{gr,ha}-\tilde{s}_{t,\kappa}R_{t}\right) + (1-\eta^\text{c}\eta^\text{d}) P_{t,\kappa}^\text{d} \notag\\
	&\leq \eta^\text{c}\left( P_{t}^\text{gr,ha}-\tilde{s}_{t,\kappa}R_{t}\right), \qquad \forall \kappa.
\end{align}
Thus, constraint (\ref{eqn_2_10}) can also be safely approximated by:
\begin{align}
	&\tilde{e}_{0} + \sum_{\kappa=1}^{d}{\eta^\text{c}\left(P_{t}^\text{gr,ha}-\tilde{s}_{t,\kappa}R_{t}\right)}\Delta d \leq \tilde{e}_t^{+}, \qquad \forall d.\label{eqn_4_9}
\end{align}

Based on the above analysis, and substitute $\sum_{\kappa=1}^{d}{\tilde{s}_{t,\kappa}}$ by $d\tilde{s}_{t,d}^\text{a}$ in the above constraints, we have that when $\tilde{p}_{t}^{+}\geq 0$, $\tilde{p}_{t}^{-}\leq 0$, and $P_{t}^\text{gr,da}$ and $P_{t}^\text{gr,ha}$ have the same sign, constraints (\ref{eqn_2_4})--(\ref{eqn_2_10}) hold if the following constraints hold: 
\begin{align}
& \eta^\text{c} \left(P_{t}^\text{gr,ha}-\tilde{s}_{t,1}^\text{a}R_{t}\right) - \tilde{p}_{t}^{+} \leq 0, \label{eqn_cc1_3}\\
& -\frac{1}{\eta^\text{d}} \left(P_{t}^\text{gr,ha}-\tilde{s}_{t,1}^\text{a}R_{t}\right) +\tilde{p}_{t}^{-} \leq 0,\label{eqn_cc2_3}\\
&\left\{\begin{array}{l}
- \eta^\text{c}d \Delta d P_{t}^\text{gr,ha} + \left( \frac{1+\eta^\text{c}\eta^\text{d}}{2\eta^\text{d}}\tilde{s}_{t,d}^\text{a} +\frac{1-\eta^\text{c}\eta^\text{d}}{2\eta^\text{d}} \right)d \Delta d R_{t}  \\
-  \tilde{e}_{0} + \tilde{e}_t^{-} \leq 0,  \quad\forall d, \quad\textit{~if~} P_{t}^\text{gr,da}\geq 0,\\
- \frac{1}{\eta^\text{d}}d \Delta d P_{t}^\text{gr,ha} + \left( \frac{1+\eta^\text{c}\eta^\text{d}}{2\eta^\text{d}}\tilde{s}_{t,d}^\text{a} +\frac{1-\eta^\text{c}\eta^\text{d}}{2\eta^\text{d}} \right)d \Delta d R_{t}  \\
-  \tilde{e}_{0} + \tilde{e}_t^{-} \leq 0,  \quad\forall d, \quad\textit{~if~} P_{t}^\text{gr,da}< 0,
\end{array}\right. \label{eqn_cc3_3}\\
& \eta^\text{c}d \Delta d P_{t}^\text{gr,ha} - \eta^\text{c}\tilde{s}_{t,d}^\text{a} d \Delta d R_{t}  + \tilde{e}_{0} -\tilde{e}_t^{+} \leq 0, \quad\forall d. \label{eqn_cc4_3}
\end{align}
Note that all the intermediate variables, i.e., $D_{t,d}, P_{t,d}^\text{c},P_{t,d}^\text{d}$ are omitted because their constraints are redundant for the true decision variables, i.e., $P_{t}^\text{gr,ha}$ and $R_{t}$. 

This proves Proposition 2. \hfill $\square$ 

We can also observe that, for the relaxed energy constraint (\ref{eqn_cc4_3}), it is slightly conservative because we relaxed the element $(1-\eta^\text{c}\eta^\text{d}) P_{t,\kappa}^\text{d}$ for each sub-hourly interval $d$. The total relaxed term $\sum_{\kappa=1}^{d}{(1-\eta^\text{c}\eta^\text{d}) P_{t,\kappa}^\text{d}} $ will increase with the increase of $d$. As a result, the new energy constraint (\ref{eqn_cc4_3}) will be more conservative when $d$ is large. Therefore, in practice, we can simply adopt constraints (\ref{eqn_cc3_3})--(\ref{eqn_cc4_3}) for $d=1800$ instead of all the possible $1\leq d\leq 1800$.

\section{The Full Formulation of the Hour-ahead Chance-constrained Program}
We give the full formulation of the hour-ahead chance-constrained program, i.e., \textbf{P3}, in this section for the convenience of future proofs:
\begin{align}
&\max \mathcal{F}= (c_t^\text{rc} + c_t^\text{rp} \overline{m}_t) R_t\Delta t -c_t^\text{e,rt} \Delta P_t^\text{gr,ha}\Delta t \notag\\
&\qquad\qquad-\sum_{\omega\in \Omega}\pi_{\omega}{c^\text{d}}E_{\omega,t}^{\text{der,d}}, \label{eqn_obj2}\\
&\text{subject to:}\notag\\
&0 \leq  R_t \leq R_t^\text{da}\label{eqn_2_0},\\
&\Delta P_t^\text{gr,ha}\geq P_t^\text{gr,ha}-P_t^\text{gr,da},\label{eqn_2_1}\\
&\Delta P_t^\text{gr,ha}\geq -P_t^\text{gr,ha}+P_t^\text{gr,da},\label{eqn_2_2}\\
&P_{t}^\text{gr,da}-s_{\omega,t}^\text{up/dn}R_{\omega,t} = \frac{P_{\omega,t}^\text{c,up/dn}}{\eta^\text{c}} + P_{\omega,t}^\text{d,up/dn}\eta^\text{d}, ~~ \forall \omega \in \Omega,\label{eqn_1_1}\\
& 0 \leq P_{\omega,t}^\text{c,up/dn} \leq p_{\omega,t}^{+}, \qquad \forall \omega \in \Omega,\label{eqn_1_3}\\
& p_{\omega,t}^{-} \leq P_{\omega,t}^\text{d,up/dn} \leq 0, \qquad \forall \omega \in \Omega,\label{eqn_1_4}\\
%& D_{\omega,t}^\text{up}, D_{\omega,t}^\text{dn} \in \{0,1\}, \qquad \forall \omega \in \Omega \label{eqn_1_7}\\
&E_{\omega,t}^\text{der,d}=-\delta t_t^\text{up}P_{\omega,t}^\text{d,up} - \delta t_t^\text{dn}P_{\omega,t}^\text{d,dn}, \qquad \forall \omega \in \Omega,\label{eqn_2_3}\\
&\text{Pr}_{\tilde{\mathbf{\xi}}^{'}\sim {\Xi}^{'}}\left\{ \eta^\text{c} \left(P_{t}^\text{gr,ha}-\tilde{s}_{t,1}^\text{a}R_{t}\right) - \tilde{p}_{t}^{+} \leq 0 \right\} \geq 1-\epsilon, \label{eqn_cc1_2}\\
&\text{Pr}_{\tilde{\mathbf{\xi}}^{'}\sim {\Xi}^{'}}\left\{ -\frac{1}{\eta^\text{d}} \left(P_{t}^\text{gr,ha}-\tilde{s}_{t,1}^\text{a}R_{t}\right) +\tilde{p}_{t}^{-} \leq 0 \right\} \geq 1-\epsilon, \label{eqn_cc2_2}\\
&\left\{\begin{array}{l}
\text{Pr}_{\tilde{\mathbf{\xi}}^{'}\sim {\Xi}^{'}}\Bigg\{ - \eta^\text{c} P_{t}^\text{gr,ha} + \left( \frac{1+\eta^\text{c}\eta^\text{d}}{2\eta^\text{d}}\tilde{s}_{t,1800}^\text{a} +\frac{1-\eta^\text{c}\eta^\text{d}}{2\eta^\text{d}} \right) R_{t} \\
-  \tilde{e}_{0} + \tilde{e}_t^{-} \leq 0 \Bigg\} \geq 1-\epsilon, \quad\textit{~if~} P_{t}^\text{gr,da}\geq 0, \\
\text{Pr}_{\tilde{\mathbf{\xi}}^{'}\sim {\Xi}^{'}}\Bigg\{ - \frac{1}{\eta^\text{d}} P_{t}^\text{gr,ha} + \left( \frac{1+\eta^\text{c}\eta^\text{d}}{2\eta^\text{d}}\tilde{s}_{t,1800}^\text{a} +\frac{1-\eta^\text{c}\eta^\text{d}}{2\eta^\text{d}} \right) R_{t} \\
-  \tilde{e}_{0} + \tilde{e}_t^{-} \leq 0 \Bigg\} \geq 1-\epsilon, \quad\textit{~if~} P_{t}^\text{gr,da}< 0, 
\end{array}\right.\label{eqn_cc3_2}\\
&\text{Pr}_{\tilde{\mathbf{\xi}}^{'}\sim {\Xi}^{'}}\left\{ \eta^\text{c} P_{t}^\text{gr,ha} - \eta^\text{c}\tilde{s}_{t,1800}^\text{a}  R_{t}  + \tilde{e}_{0} -\tilde{e}_t^{+} \leq 0\right\} \geq 1-\epsilon. \label{eqn_cc4_2}
\end{align}

\section{Reformulation of \textbf{P3}}
In this section, we reformulate \textbf{P3} into its standard form.

Based on equations (\ref{eqn_cc1_3})--(\ref{eqn_cc4_3}), we define $\tilde{\mathbf{a}}_j$ and $\tilde{b}_j$ as:
\begin{align}
&\tilde{\mathbf{a}}_j =     \notag\\
&\left\{ 
\begin{array}{l}
         {\left[-\eta^\text{c}\tilde{s}_{t,1}^\text{a}, ~ \eta^\text{c}\right]}^{\intercal},
         \qquad j=1 \\ 
         {\left[\frac{1}{\eta^\text{d}}\tilde{s}_{t,1}^\text{a}, ~-\frac{1}{\eta^\text{d}} \right]}^{\intercal}, \qquad j=2 \\
         {\left[ \frac{1+\eta^\text{c}\eta^\text{d}}{2\eta^\text{d}}\tilde{s}_{t,d}^\text{a} +\frac{1-\eta^\text{c}\eta^\text{d}}{2\eta^\text{d}},~ - \eta^\text{c}\right]}^{\intercal}, ~P_{t}^\text{gr,da}\geq 0, ~ j=3\\
         {\left[ \frac{1+\eta^\text{c}\eta^\text{d}}{2\eta^\text{d}}\tilde{s}_{t,d}^\text{a} +\frac{1-\eta^\text{c}\eta^\text{d}}{2\eta^\text{d}},~ - \frac{1}{\eta^\text{d}}\right]}^{\intercal}, ~ P_{t}^\text{gr,da}< 0,~~ j=3\\
         {\left[ - \eta^\text{c}\tilde{s}_{t,d}^\text{a}, ~\eta^\text{c}\right]}^{\intercal}, \qquad j=4 \\
\end{array},\right.\label{aj}\\
&\tilde{b}_j =     \left\{ 
\begin{array}{l}
	- \tilde{p}_{t}^{+}, \qquad j=1 \\ 
	\tilde{p}_{t}^{-}, \qquad j=2 \\
	-  \tilde{e}_{0} + \tilde{e}_t^{-}, \qquad j=3 \\
	\tilde{e}_{0} -\tilde{e}_t^{+}, \qquad j=4 \\
\end{array},\right.\label{bj}
\end{align}
which can reformulate the constraints (\ref{eqn_cc1_3})--(\ref{eqn_cc4_3}) by $\tilde{\mathbf{d}}_j^{\intercal}\hat{\mathbf{X}}\leq 0$, in which $\tilde{\mathbf{d}}_j \doteq [\tilde{\mathbf{a}}_j^{\intercal} ~ \tilde{b}_j]^{\intercal}$;  $\hat{\mathbf{X}}\doteq[\mathbf{X}^{\intercal}~ 1]^{\intercal}$; $\mathbf{X}\doteq[R_t,  P_t^\text{gr,ha} ]^{\intercal}$.

As a result, \textbf{P3} can be reformulated as:
\begin{align}
\textbf{P3\_st:}\quad &\max \quad \text{(\ref{eqn_obj2})} \qquad \text{s.t.:\quad (\ref{eqn_2_0})-(\ref{eqn_2_3}), and }\notag\\
&\text{Pr}_{\tilde{\mathbf{\xi}}^{'}\sim {\Xi}^{'}}\left\{ \tilde{\mathbf{d}}_j^{\intercal}\hat{\mathbf{X}}\leq 0 \right\} \geq 1-\epsilon, ~ \forall j=1,2,3,4,
\end{align}
This is the standard formulation for the linear chance-constrained program \textbf{P3}. 

Note that, in \textbf{P3\_st}, we use the new uncertain parameters $\tilde{\mathbf{\xi}}^{'}$ instead of $\tilde{\mathbf{\xi}}$. This process leads to two benefits: 1) First, the scale of each vector of uncertain parameters, i.e., $\tilde{\mathbf{d}}_j $, is significantly decreased (by contrast, if we use $\tilde{s}_{t,d}$, the scales of uncertain parameters in (\ref{eqn_cc3_2})--(\ref{eqn_cc4_2}) increase with $d$); 2) Second, all the entities in $\tilde{\mathbf{d}}_j$ are independent with each other (by contrast, $\tilde{s}_{t,d}$ in adjacent intervals are dependent). As a result, it is much easier to calculate the covariance of $\tilde{\mathbf{d}}_j$ (see the following section). 

\section{Calculation of the Mean and Covariance in \textbf{P4} }
The mean and covariance of $\tilde{\mathbf{d}}_j$, i.e., $\overline{\mathbf{d}}_j$ and $\mathbf{\Gamma}_j$, can be calculated as follows:
\begin{align}
&\overline{\mathbf{d}}_j =     \notag\\
&\left\{ 
\begin{array}{l}
{\left[-\eta^\text{c}\overline{s}_{t,1}^\text{a}, ~ \eta^\text{c},~ - \overline{p}_{t}^{+}\right]}^{\intercal},
\qquad j=1 \\ 
{\left[\frac{1}{\eta^\text{d}}\overline{s}_{t,1}^\text{a}, ~-\frac{1}{\eta^\text{d}}, ~\overline{p}_{t}^{-} \right]}^{\intercal}, \qquad j=2 \\
{\left[ \left(\frac{1+\eta^\text{c}\eta^\text{d}}{2\eta^\text{d}}\overline{s}_{t,1800}^\text{a} +\frac{1-\eta^\text{c}\eta^\text{d}}{2\eta^\text{d}}\right),~ - \eta^\text{c},~-  \overline{e}_{0} + \overline{e}_t^{-}\right]}^{\intercal},  \\
\qquad\qquad\qquad\qquad\qquad \qquad ~\text{if}~P_{t}^\text{gr,da}\geq 0,  j =3\\
{\left[ \left(\frac{1+\eta^\text{c}\eta^\text{d}}{2\eta^\text{d}}\overline{s}_{t,1800}^\text{a} +\frac{1-\eta^\text{c}\eta^\text{d}}{2\eta^\text{d}}\right),~ -\frac{1}{\eta^\text{d}},~-  \overline{e}_{0} + \overline{e}_t^{-}\right]}^{\intercal},  \\
\qquad\qquad\qquad\qquad\qquad \qquad ~\text{if}~P_{t}^\text{gr,da}< 0,  j =3\\
{\left[ - \eta^\text{c}\overline{s}_{t,1800}^\text{a}, ~\eta^\text{c},~ \overline{e}_{0} -\overline{e}_t^{+}\right]}^{\intercal}, \qquad j=4 \\
\end{array},\right.\label{dj}\\
&\mathbf{\Gamma}_j =     \notag\\
&\left\{ 
\begin{array}{l}
{\text{diag}\left(\eta^\text{c}\text{Var}(\tilde{s}_{t,1}^\text{a}), ~ 0,~  \text{Var}(\tilde{p}_{t}^{+})\right)},
\qquad j=1 \\ 
{\text{diag}\left(\frac{1}{\eta^\text{d}}\text{Var}(\tilde{s}_{t,1}^\text{a}), ~0, ~\text{Var}(\tilde{p}_{t}^{-}) \right)}, \qquad j=2 \\
{\text{diag}\left(\frac{1+\eta^\text{c}\eta^\text{d}}{2\eta^\text{d}}\text{Var}(\tilde{s}_{t,1800}^\text{a}),~ 0,~ \text{Var}(\tilde{e}_{0}) + \text{Var}(\tilde{e}_t^{-})\right)},  \\
\qquad\qquad\qquad\qquad \qquad\qquad\qquad\qquad  j =3\\
{\text{diag}\left( \eta^\text{c}\text{Var}(\tilde{s}_{t,1800}^\text{a}), ~0,~ \text{Var}(\tilde{e}_{0}) +\text{Var}(\tilde{e}_t^{+})\right)},  ~~ j=4 \\
\end{array},\right.\label{gammaj}
\end{align}
where, we use over-line symbols, e.g., $\overline{s}_{t,1}^\text{a}$, and Var($\cdot$), e.g., $\text{Var}(\tilde{s}_{t,1}^\text{a})$, to denote the mean and variance of the uncertain parameters in $\tilde{\mathbf{\xi}}^{'}$, which can be estimated learning historical data; diag($\cdot$) denotes a diagonal matrix. Because the variance of a uncertain parameter is nonnegative, we have $\mathbf{\Gamma}_j\succeq 0$. 

\section{Distributions of the Regulation Signals}
The regulation signals (RegA-type) from 1st August, 2015, to 31th July, 2016 (a total of 366 days), of the PJM market\cite{PJM_pilong2013pjm} are used to analyze their distributions. 

\begin{figure}
	\centering
%	\vspace{-4mm}
	\begin{minipage}[t]{1\columnwidth}
		\centering
		\epsfig{file = 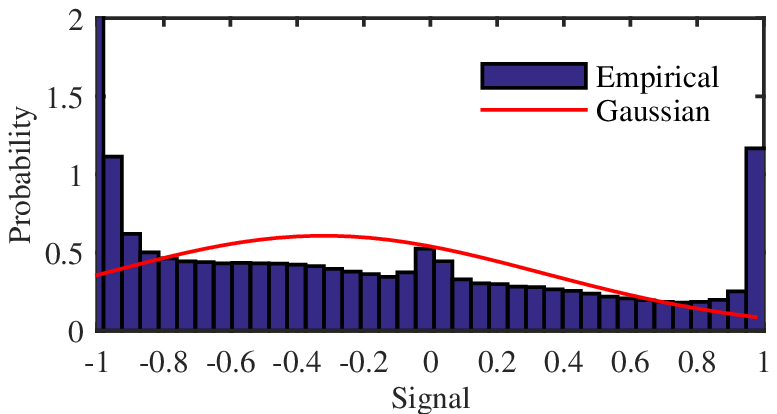, width = 6.75cm}
		\vspace{-3mm}
		\caption{Distribution of $\tilde{s}_{t,1}^\text{a}$.}
		\label{fig_second}
		\vspace{1mm}
	\end{minipage}
	\begin{minipage}[t]{1\columnwidth}
		\centering
		\epsfig{file = 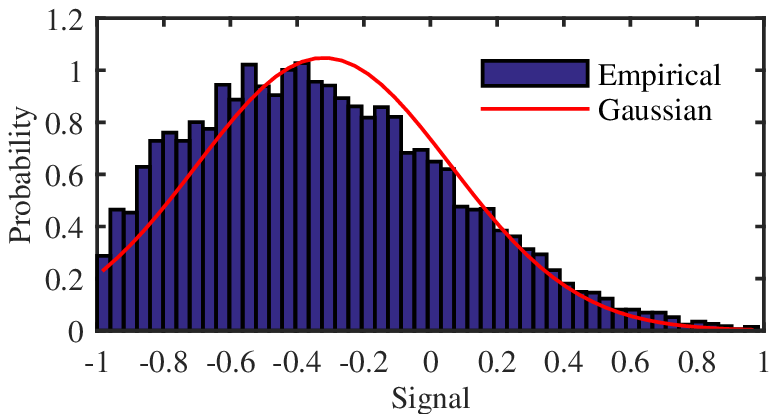, width = 6.75cm}
		\vspace{-3mm}
		\caption{Distribution of $\tilde{s}_{t,1800}^\text{a}$.}
		\label{fig_hour}
	\end{minipage}
	\vspace{-4mm}
\end{figure}

The empirical distributions and its Gaussian regressions of the regulation signals $\tilde{s}_{t,1}^\text{a}$ and $\tilde{s}_{t,1800}^\text{a}$ are respectively illustrated in Fig.~\ref{fig_second} and Fig.~\ref{fig_hour}. It is obvious that the regulation signals in every two seconds, i.e., $\tilde{s}_{t,1}^\text{a}$, do not follow a specific tractable distribution; by contrast, the distribution of $\tilde{s}_{t,1800}^\text{a}$ is approximately Gaussian.

\section{Proof of Proposition 3}
In this section, we briefly prove the proposition 3.

The data-driven distributionally robust chance-constrained program is:
\begin{align}
&\textbf{P4:}~\max \quad \text{(\ref{eqn_obj2})} \qquad \text{s.t.:\quad (\ref{eqn_2_0})-(\ref{eqn_2_3}), and }\notag\\
&\inf_{\mathcal{D}_j \in (\overline{\mathbf{d}}_j, \mathbf{\Gamma}_j)}\text{Pr}_{\tilde{\mathbf{d}}_j\sim \mathcal{D}_j}\left\{ \tilde{\mathbf{d}}_j^{\intercal}\hat{\mathbf{X}}\leq 0 \right\} \geq 1-\epsilon,  ~~\forall j =1, 2, \label{eqn_CC3}\\
&\inf_{ I_{\phi}(\mathcal{D}_j, \mathcal{D}_0)\leq \rho}\text{Pr}_{\tilde{\mathbf{d}}_j\sim \mathcal{D}_j}\left\{ \tilde{\mathbf{d}}_j^{\intercal}\hat{\mathbf{X}}  \leq 0 \right\} \geq 1-\epsilon, ~\forall  j=3,4. \label{eqn_CC4}
\end{align}

\textbf{Proposition 3} \textit{For any $\epsilon \in (0,0.5]$, the chance constraints (\ref{eqn_CC3})-(\ref{eqn_CC4}) hold, if the following constraints hold} 
\begin{align}
&\sqrt{{(1-\epsilon)}/{\epsilon}} \| \mathbf{\Gamma}_j^{\frac{1}{2}} \hat{\mathbf{X}}\|_2  + \overline{\mathbf{d}}_j^{\intercal} \hat{\mathbf{X}} \leq 0, \qquad \forall j=1,2,\label{eqn_CC3_2}\\
&\Psi_G^{-1}(1-\epsilon^{'}
) \| \mathbf{\Gamma}_j^{\frac{1}{2}} \hat{\mathbf{X}}\|_2  + \overline{\mathbf{d}}_j^{\intercal} \hat{\mathbf{X}} \leq 0, \qquad \forall j=3,4,\label{eqn_CC4_2}
\end{align}
\textit{where, $\Psi_G(\epsilon)=(1/\sqrt{2\pi})\int_{-\infty}^{\epsilon}\text{exp}(\frac{t^2}{2})\text{d}t$; $\epsilon^{'}=\epsilon-\frac{\sqrt{\rho^2 +4\rho(\epsilon^2-\epsilon)}-(1-2\epsilon)\rho}{2\rho+2}$ is an adjusted tolerance based on $\chi^2$-divergence. }

\textit{Proof:} First, based on the result proposed by Calafiore and Ghaoui \cite{DCC_Calafiore2006} (Theorem 3.1), for any $\epsilon \in (0,1)$, the chance constraint (\ref{eqn_CC3}) holds, if  constraint (\ref{eqn_CC3_2}) holds.

To solve the $\phi$-divergence based data-driven chance-constraint (\ref{eqn_CC4}), we introduce another proposition based on the findings of Jiang \& Guan \cite{DCC_Data_Driven_Jiang2016} (we omitted the proof here for brevity, interested readers can refer to Corollary 1 and Proposition 2 in \cite{DCC_Data_Driven_Jiang2016}):

\textbf{Proposition 4} \textit{(Jiang \& Guan\cite{DCC_Data_Driven_Jiang2016}, Corollary 1 and Proposition 2): The distributionally robust chance constraint}
	\begin{align}
	&\inf_{ I_{\phi}(\mathcal{D}_j, \mathcal{D}_0)\leq \rho}\text{Pr}_{\tilde{\mathbf{d}}_j\sim \mathcal{D}_j}\left\{ \tilde{\mathbf{d}}_j^{\intercal}\hat{\mathbf{X}}  \leq 0 \right\} \geq 1-\epsilon, ~~\forall  j, \notag
	\end{align}
	\textit{where, $I_{\phi}(\mathcal{D}_j, \mathcal{D}_0)$ is defined by the $\chi^2$-divergence,
	holds if the classical chance constraint}
	\begin{align}
	&\text{Pr}_{\tilde{\mathbf{d}}_j\sim \mathcal{D}_0}\left\{ \tilde{\mathbf{d}}_j^{\intercal}\hat{\mathbf{X}}  \leq 0 \right\} \geq 1-\epsilon^{'}, ~~\forall  j, \notag
	\end{align}
	\textit{holds, where}
	\begin{align}
	&\epsilon^{'} = \epsilon-\frac{\sqrt{\rho^2 +4\rho(\epsilon^2-\epsilon)}-(1-2\epsilon)\rho}{2\rho+2}.\notag
	\end{align}

By Proposition 4, we can conservatively approximate the chance constraint (\ref{eqn_CC4}) by the following constraint:
\begin{align}
&\text{Pr}_{\tilde{\mathbf{d}}_j\sim \mathcal{D}_0}\left\{ \tilde{\mathbf{d}}_j^{\intercal}\hat{\mathbf{X}}  \leq 0 \right\} \geq 1-\epsilon^{'}, ~~\forall  j=3,4, \label{eqn_DCC_gassian}
\end{align}
in which, $\mathcal{D}_0$ is Gaussian with mean $\overline{\mathbf{d}}_j$ and variance $\mathbf{\Gamma}_j$.

Then, based on the result on Gaussian chance-constrained program proposed by Pr{\'e}kopa \cite{DCC_prekopa1995stochastic} (Theorem 10.4.1), for any $\epsilon\in (0,0.5]$, the chance constraint (\ref{eqn_DCC_gassian})
holds, if the following constraint holds:
\begin{align}
&&\Psi_G^{-1}(1-\epsilon^{'}
) \| \mathbf{\Gamma}_j^{\frac{1}{2}} \hat{\mathbf{X}}\|_2  + \overline{\mathbf{d}}_j^{\intercal} \hat{\mathbf{X}} \leq 0,  ~~\forall j =3, 4,\notag
\end{align}
which is exactly the constraint (\ref{eqn_CC4_2}).

This proves Proposition 3. \hfill $\square$

\bibliographystyle{ieeetr}
\bibliography{reference}